\documentclass[11pt,a4paper]{article}

\usepackage[T1]{fontenc}
\usepackage[utf8]{inputenc}
\usepackage{upgreek}
\usepackage{hyperref}
\usepackage{graphicx}
\usepackage{multicol}
\usepackage{lipsum}
\usepackage{geometry}
\usepackage{amsmath, amssymb, amsthm}

\usepackage[capitalise,nameinlink,sort&compress]{cleveref}
\crefname{equation}{}{}
\crefname{appendix}{Appendix}{appendix}
\creflabelformat{equation}{#2#1#3}
\crefmultiformat{equation}{#2#1#3}{, #2#1#3}{, #2#1#3}{, #2#1#3}
\renewcommand\eqref[1]{(\cref{#1})}

\usepackage{xcolor}\usepackage{xspace}
\newcommand{\cor}[1]{\textcolor{black}{#1}\xspace}

\title{The illusion of illusions:\\ There are no optical corrections in the Parthenon}
\author{Alain Goriely\\
\small Mathematical Institute, University of Oxford, UK\\
\small \texttt{alain.goriely@maths.ox.ac.uk}}
\date{}

\begin{document}
\maketitle

\begin{abstract}
One of the oldest and most enduring myths in human history is the belief that the Parthenon was cleverly designed with various curved structures and sizes in order to correct optical illusions, and therefore appear straight and regular.  The myth has its origin in  the writings of Vitruvius more than 2,000 years ago and was renewed in the nineteenth century when curved profiles were carefully measured. At least twelve different \textit{optical corrections}  have been proposed, all following the same basic principle. The myth is still widely acknowledged as an obvious truth despite a complete absence of historical or scientific evidence.  This paper analyzes   these corrections scientifically and demonstrates that the illusions they are supposed to correct are either non-existent or  so small as to  be imperceptible. 
\end{abstract}

 \section{Introduction}
\subsection{The Parthenon}
The Parthenon is a classical Greek temple situated on the Acropolis hill overlooking Athens, Greece. Built   during the Golden Age of Athens under Pericles, it was completed in 15 years, between 447 and 432 BCE.  Dedicated primarily to Athena Parthenos, the city’s patron goddess, it was designed by the architects Ictinus and Callicrates under the  supervision of Phidias, the sculptor.  The rectangular structure measures approximately $a=30.9\,\text{m}$ in width aligned approximately along the north--south axis, and $b=69.5\,\text{m}$ in length along the east--west axis.  In terms of architecture, it is mostly classified as  Doric architecture (with its typical simple columns), but also contains Ionic elements. Over time, the Parthenon has served various purposes, functioning as a temple, treasury, church, mosque, and unfortunately, as a gunpowder depot, which led to extensive damage when Venetian forces bombarded the Acropolis on 26 September 1687,  causing  the collapse of the roof, parts of the interior walls, and numerous columns. Further damage occurred when Thomas Bruce, the 7th Earl of Elgin, removed significant portions of the Parthenon sculptures between 1801 and 1812, which can now be found in the British Museum. Yet, it  still stands and  is widely regarded as the highest achievement of Greek art and architecture and an iconic monument of classical antiquity.

\begin{figure}[ht!]
 \centering
\includegraphics[width=0.8\linewidth]{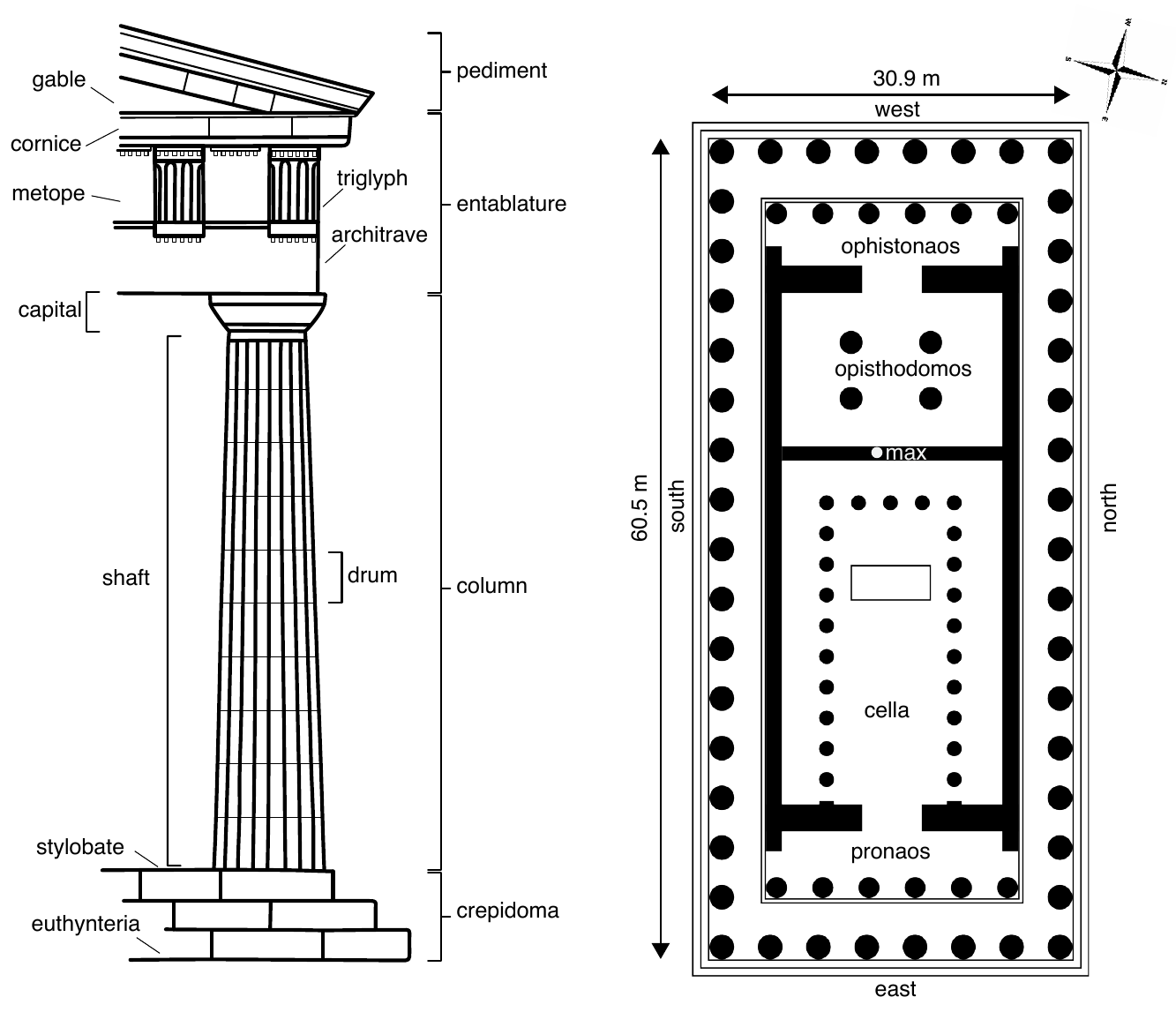}
 \caption{Left: The classical architectural elements of a Doric temple. Right: Plan of the Parthenon. }
 \label{doric-temple}
\end{figure}

Greek temples like the Parthenon were built  with distinct architectural elements on a stepped base, the \textit{crepidoma} (\cref{doric-temple}). This base  typically has  three steps that elevate the temple above the surrounding ground. The  \textit{euthynteria} sits right on the foundation. The \textit{stylobate},  of particular interest for our discussion, is the top platform (of width $a$ and length $b$) on which  vertical elements  stand. It  supports the \textit{peristyle},  the outer ring of columns, and the \textit{cella}, which is the central chamber containing the cult statue, the heart of the building.

 The peristyle  provides the main vertical support for the top horizontal superstructure, the \textit{entablature}. The entablature itself is classically divided into three  parts: the \textit{architrave}, a beam that rests directly on the \textit{capitals} of the columns; a horizontal band with either \textit{triglyphs} and \textit{metopes} in the Doric order or a \textit{frieze} decorated with relief sculptures in the Ionic order; and the \textit{cornice}, the projecting uppermost section that frames and shelters the structure beneath. These elements support the \textit{pediment}, the triangular gable at each end of the temple, which itself is often filled with decoration. The Parthenon preserves the Doric triglyph–metope rhythm on its exterior but departs from it internally through the addition of a continuous Ionic frieze. 

\subsection{The refinements}

One of the most striking features of the Parthenon is that despite the appearance that it is built from straight vertical and horizontal periodic elements, most features slightly deviate from their perceived regularity. It is often said, erroneously, that no straight line can be found in the Parthenon.  Dozens of these \textit{architectural refinements} have been identified over the years, including the upward curvature of the stylobate and entablature, the slight bulging of column shafts (\textit{entasis}),  the inward inclination of columns, and the contraction of column spacing at the corners, among many others. These refinements are sometimes called \textit{optical refinements} but this term will not be used here as it wrongly attributes a purpose for these variations and therefore naturally propagates the myth. For instance,  the underground foundation of the Parthenon is also curved, hence refined despite being invisible. It is now accepted is that many of these refinements were executed with precision and occur consistently across structures,  indicating that they were not construction flaws but deliberate variations from rectilinear and regular forms. Our modern understanding of civil engineering further indicates that they are not necessary for structural integrity and a lingering question for more than two millennia has been to understand the motivation behind these architectural refinements.

\subsection{The myths}

\begin{figure}[ht!]
 \centering
\includegraphics[width=\linewidth]{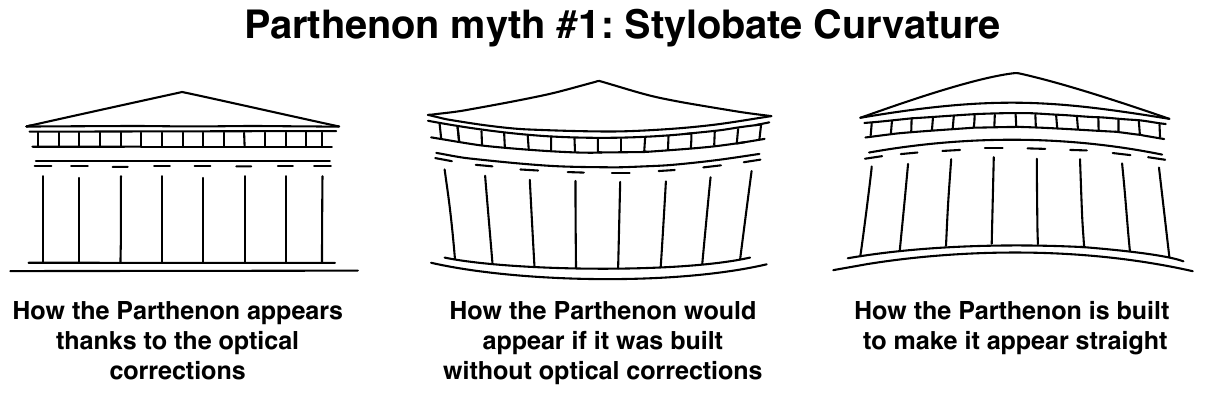}
 \caption{The myth of the Parthenon follows a typical reasoning: In the absence of refinements, the Parthenon would look like it is sagging as in the central figure (grossly exaggerated for the purpose of the argument). But, the ancient Greeks, in their great wisdom, built the Parthenon with upward curvature (right) so that it would appear straight (left). All three statements are incorrect. Figure adapted from Fletcher's \cite{fletcher1905history}. }
 \label{fig1}
\end{figure}

The enduring myth is that architectural refinements were \textit{optical corrections},  introduced deliberately with the explicit intention of counteracting perceived visual illusions and restoring regularity or rectilinearity. According to the myth, buildings were designed with curves so that they look straight and with irregular spacing so that they appear regular. The paradigm for this idea, shown in \cref{fig1}, is the   upward curvature of the stylobate and the entablature, explained by stating that a building with straight horizontal features would appear to sag. Hence, it is claimed that an upward curvature corrects this alleged illusion, and the building appears straight. As will be shown, this explanation, like all such contrived arguments, collapses at the very first scrutiny. Yet, these myths are so universal that they have been repeated in specialized books about the Parthenon, encyclopedias, documentaries, and, ad nauseam, over the planetary dumping ground of human knowledge, propaganda, truths, lies, and myths, called the internet. A simple measure of its widespread acceptance is found by asking any Large Language Model (such as ChatGPT, Gemini, Claude, all interrogated on 9 September 2025) the simple question: what are optical refinements? Invariably, these systems will confidently give the  correction theory as the reason for the refinements and cite  the Parthenon as its prime example. This is not surprising, as large language models are trained on vast amounts of text that reflect what has been written and repeated over centuries. Since the optical illusion explanation for the Parthenon has been taught in art history and architecture courses, repeated across textbooks and encyclopedias, cited in museum descriptions and tour guides, documentaries, travel books, social media, and  assumed as fact in countless articles and books, any AI system will reproduce this explanation for the simple reason that it pervades its training data.\\

The main myths associated with the Parthenon are:
\begin{itemize}
    \item \textbf{Myth \#1: Stylobate Curvature.}
Horizontal elements such as the stylobate, entablature, and pediment are curved, preventing them from appearing to sag in the middle. 
\item \textbf{Myth \#2: Column Entasis.}
Columns are not perfectly straight but have a slight bulge in the middle, preventing them from appearing concave or pinched. 
\item \textbf{Myth \#3: Corner Column Thickening.}
Corner columns are made slightly thicker  than the intermediate ones, offsetting the stronger light around the edges that would make them seem thinner.
\item \textbf{Myth \#4: Proportional Adjustments.} The proportional enlargement of elements placed at greater height, such as texts or human heads, to make them appear of uniform size with elements placed at lower height.
\item \textbf{Minor Myths.}
Other myths include the slight inward tilt of  
columns, supposedly creating the illusion of perspective and preventing the building from appearing to bulge outward at the top, and  intercolumnation: columns at the corners are often closer together, supposedly counteracting the visual effect of a wider space between them when viewed from a distance. 
\end{itemize}

\subsection{A brief history}
Ancient Greeks were aware of the existence of visual illusions as some were already described by Plato and Aristotle \cite{johannsen1971early}. However, none of the illusions that have been cited as possibly occurring in the Parthenon are to be found in the ancient literature. We have to wait 400 years after the Parthenon's construction for the myth to take shape. Indeed,  the earliest and main historical source is Vitruvius, an important but notoriously opaque writer. In \textit{De architectura} (first century BCE), the only architectural text we have from antiquity, he states that certain refinements, such as   the slight curvatures in stylobates (myth \#1),  column entasis (myth \#2), larger corner columns (myth \#3),  and scale corrections (myth \#4), were made to counteract optical illusions, specifically the appearance of sagging, concavity in long horizontal lines, or changes in size. However, while Vitruvius may have had access to some of the earlier Greek writings, as argued by John James Coulton \cite[p. 109]{coulton1982ancient}, he provided no evidence that these refinements were intended by the Greeks as optical corrections. 
 
During the Renaissance, the rediscovery of Vitruvius gave renewed authority to the idea of optical corrections. Humanists and architects such as Alberti, Serlio, and Palladio \cite{Alberti_DeReAedificatoria_1485,Serlio_TutteLOpere_1537,Palladio_QuattroLibri_1570} drew on Vitruvian texts mostly to discuss scale corrections \cite{Frangenberg1993}.
But, it is really in the nineteenth century that the myth took hold following  a new phase of precision measurements. The most authoritative and enduring formulation of the optical-correction theory was advanced by Francis Cranmer Penrose in his majestic \textit{An Investigation of the Principles of Athenian Architecture} (1851, revised 1888) \cite{Penrose1851}, which became the foundational reference for subsequent interpretations of the Parthenon’s design. The book itself, published by the \textit{Society of Dilettanti}, is truly monumental, measuring 570 $\times$ 410 mm with color engravings, fold-out plates, and detailed  drawings. Penrose first demonstrated, with unprecedented accuracy, the presence of deliberate refinements. Then, he argued that such refinements could only have been introduced as optical corrections. Penrose’s interpretation  was further expanded by John Pennethorne in \textit{The Geometry and Optics of Ancient Architecture} (1878)\cite{Pennethorne1878} and  repeated and systematically popularized through highly influential textbooks such as Banister Fletcher’s \textit{History of Architecture on the Comparative Method} (1896 and later editions)\cite{Fletcher1905}, Auguste Choisy  \textit{Histoire de l’architecture} (1899) \cite{Choisy1899}, and William Bell Dinsmoor’s \textit{Architecture of Ancient Greece} (1902 onward)\cite{Dinsmoor1973}. These authors systematized  and expanded the Vitruvius–Penrose explanation for a modern  readership with attractive figures, ensuring that the myth remained a standard part of architectural-history teaching well into the twentieth century.

Despite its public prevalence, the myth was never universally accepted among experts. Already in the seventeenth century, Claude Perrault,  in his \textit{Ordonnance des cinq espèces de colonnes} \cite[p.~108]{Perrault1683}  argued  that there is "\textit{no optical reason to change the proportions}". Furthermore, in his own massively annotated translation of Vitruvius, Perrault is  critical of many statements  \cite{Perrault1673}, implying that Vitruvius' explanations lack empirical or rational basis and that refinements  are  mainly due to aesthetics.
Later, by the end of the nineteenth century, Auguste Choisy in his influential \textit{Histoire de l’architecture} (1899) enthusiastically adopts  the optical correction theory and further states that both Egyptians and Greeks were fully aware and played with these illusions. However, in his contradictory conclusion he notes that the curvatures can actually be observed and that maybe their purpose was not to give a rectilinear appearance but, on the contrary,  ``\textit{to escape the vulgar aspect of buildings with straight lines''}\cite[Vol. 1, p. 324]{Choisy1899}. But the most salient rebuttal came in 1912 from John  Goodyear  \cite{goodyear1912greek} who coined \textit{architectural refinement}  and dedicated an entire book on the subject, proposing instead that these refinements may have been symbolic, structural, or purely aesthetic. 

Today, while the myth is put in doubt by some experts, it is taken for granted by the larger public. Indeed, most historians acknowledge that not all  refinements may   be for optical corrections. Yet, probably due to peer pressure, even the fiercest critics believe in some of them. For instance,  Maurice H. Pirenne, known for his work in visual perception and art, like Goodyear dismissed the idea of optical correction but  acknowledged that change in column diameters ``\textit{was to avoid their looking too thin against the bright sky''} \cite[p.~150]{Pirenne1970}. Mannolis Korres, a major Greek restoration architect and current head of the Acropolis Restoration Service, admits that the theory of optical corrections \textit{``does not hold water in most cases"} but does believe that it plays a role in the triangular pediment \cite{Korres1998}. Similarly,  Elizabeth Rankin  rejects most claims of optical corrections with clear and pointed criticism. She nevertheless conceded that the thickening of the corner columns may constitute a credible case \cite{Rankin1986}. \cor{Yet, none of these authors has presented direct geometric or scientific evidence for the presence or absence of illusions.}

The thesis of this paper  is that  there is insufficient historical evidence and strictly no scientific evidence for the extraordinary claim that the refinements of the Parthenon were either implemented as optical corrections or achieved the purpose of optical corrections. By extension, the same conclusion applies to almost all classic and modern buildings. Apart from a handful or modern buildings where architects have deliberately played with visual illusions no evidence supports the existence or necessity of optical corrections. The remainder of this paper systematically examines these refinements and presents critical and  scientific arguments against the notion that they constitute optical corrections. Alternative theories are also briefly discussed, along with the evidence that would be required to demonstrate the existence of optical corrections.

\section{Myth \#1: Curved Horizontal Elements}

\subsection{The refinement}
Following the careful measurements of the stylobate  in the first part of the nineteenth century by Joseph Hoffer in 1838 \cite{Hoffer1838}, John Pennethorne (with notes circulating in 1844 \cite{Pennethorne1844} but remained unpublished until 1878 \cite{Pennethorne1878}), and independently by Francis Penrose in 1851 \cite{Penrose1851},  it became increasingly clear that the elements found in the Parthenon, like many other Doric temples, have a slight deliberate upward curvature, a view that is now accepted by all experts \cite{haselberger1999appearance}. It is found not only on the stylobate, but also in the euthynteria that is mostly below ground and translated to the entablature by the columns. 

There has been a controversy over the years on the question of the exact nature of this curve, whether it is a segment of a large circle, of an ellipse, a parabola, or an inverted catenary \cite{stevens1943curve}.  The problem is both trivial and intellectually compelling. It is trivial in the sense that the precise mathematical form of the curve is of little consequence: given the minute degree of curvature, it can be represented with equal accuracy by several functions. A parabolic expression is convenient, and higher-order corrections to this nearly perfect fit are neither necessary for accuracy nor significant for interpretation. Indeed,  the four sides of the stylobate are exceptionally well approximated by a circular arc or, equivalently, by a parabola \cite{Georgopoulos2012}. Furthermore, it is important to point out that the first recorded systematic description of conics such as parabolas and ellipses is attributed to Menaechmus (ca. 380--320 BCE) \cite{Keppens2021}, almost a century after the construction of the Parthenon, a fact already advanced by the  Greek mathematician Constantin Carathéodory in 1937 \cite{Georgopoulos2012}. What is genuinely intriguing is the question of how the ancient Greeks conceived and executed these subtle curves \cite{coulton1993coping}. The simplest method is to pull a string between the two ends to create a catenary. It would have a natural sag controlled by end tension and the measurements could be reflected with respect to a level line to obtain the required corrections. A parabola could be created by using small intermediate steps in the spirit of Vitruvius \cite{Stevens1934}. A more complicated method found in use for the columns of Didyma is to project an elliptical arc \cite{haselberger1985construction}.

\begin{figure}[ht!]
 \centering
\includegraphics[width=\linewidth]{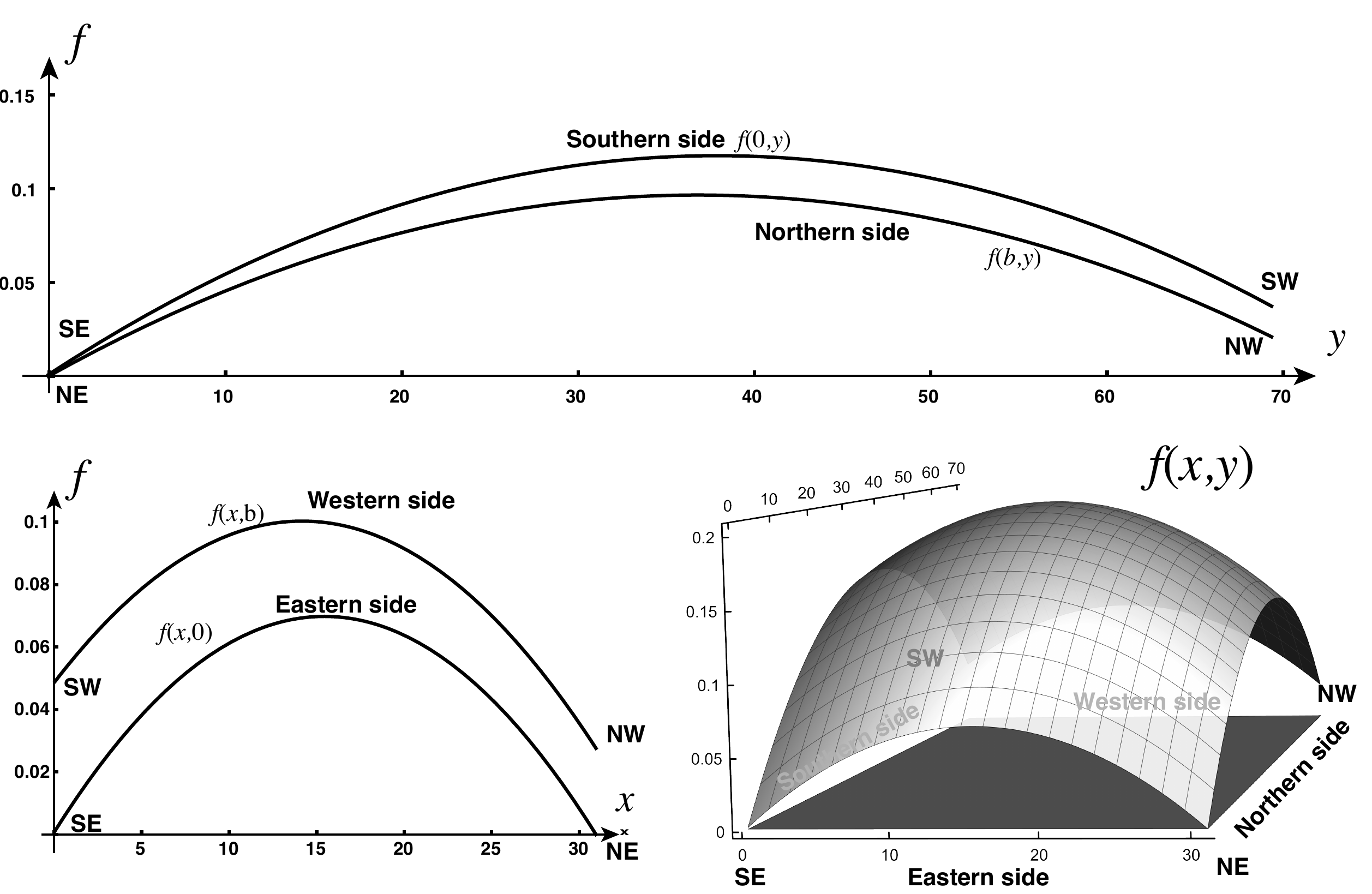}
 \caption{The shape of the four sides of the stylobates and their interpolation in 3D. Note  the huge difference in vertical scale compared to the horizontal scales and also that the four corners are not at the same level with a small overall increasing slope in the western direction. }
 \label{Fig-stylo}
\end{figure}

It is interesting to obtain an estimate of the entire stylobate profile. Using the data from Georgopoulos \cite{Georgopoulos2012}, we can easily obtain an interpolation of the entire stylobate that matches the four recorded sides  by using Coons' patch method \cite{coons_1967_coons_patch} to obtain:
\begin{align}\label{fxy}
\nonumber    f(x,y)=10^{-3}\times\Big(&1.16078+8.87918 x+8.00577 y-0.288568 x^2-0.105338
   y^2 \\
\nonumber  & -0.0523071 x y+0.000459584 x^2 y+0.000413477 x y^2\Big),\\
  &\qquad\qquad\qquad\qquad\qquad\qquad\qquad\qquad\qquad\qquad x\in[0,a],\ y\in[0,b],
\end{align}
where, as before, $a=30.9$ m, $b=69.5$ m. Here, $x=y=0$ denotes the south-east corner of the stylobate, $x$ and $y$ measuring the distance towards the northern and western directions, respectively.  In particular, the four sides are given by parabolas: $f(x,0)$ fits the east facade, $f(x,b)$, the west facade, $f(0,y)$ the south facade, and $f(a,y)$ the north facade as shown in \cref{Fig-stylo}.

By solving the coupled system $\partial_xf(x,y)=\partial_y f(x,y)=0$ for $(x,y)$, we obtain the maximum of this function, the \textit{crown} of the platform, located at coordinates $(x=13.9,y=37)$. This maximum is indicated by a white dot in the wall separating the cella and the opisthodomos in \cref{doric-temple}. It is unclear whether the location of the maximum precisely on this wall is intentional or coincidental. 

Finally, it is important to note that the upward curvature of horizontal elements found in ancient temples should not be confused with  the engineering practice of \textit{cambering}  which is the structural use of  upward curvature in  beams, girders, and slabs so that they become horizontal under expected deformation due to  loads \cite{BarkerPuckett2013}. For instance, the Golden Gate Bridge in San Francisco was built with an upward camber amounting to about 1.7 meter over a main span of 1,280 meters. Despite that, under loads it sags in the middle by about 2 meters \cite{game2016full}. Clearly, no cambering is taking place in the Parthenon's stylobate.

\subsection{The myth}
Now that we have established the existence of this curve, we come to the principal myth: the belief that  horizontal elements appear to sag and hence that this upward curvature is an optical correction so that they would appear perfectly straight as shown in \cref{fig1}. It was Vitruvius who first offered this explanation in Book 3.4.5 of  \textit{De Architectura}. The Gwilt translation of 1874 \cite{Gwilt1874} reads:
\textit{``The stylobata should be so adjusted, that, by means of small steps of stools (scamilli impares), it
may be highest in the middle. For if it be set out level, it will have the appearance of having sunk in the centre."
}\footnote{Art historians will know that there are no two words in ancient writings that have created more confusion than the infamous \textit{scamilli impares} \cite{campbell1980scamilli} (roughly translated as unequal or irregular steps) and entire books have been devoted to the subject as early as 1612 \cite{Baldi1612}. The problem is compounded by the fact that Vitruvius refers to a figure in order to explain the method. Unfortunately, the figure has been lost and the exact meaning of scamilli impares was left forever open to interpretation.} 

It is not clear when this myth reappeared in modern times as it was  only acknowledged indirectly by Penrose citing Vitruvius that \textit{``If the
line of the stylobate were perfectly horizontal, it would appear like the bed of a channel"} and concluding that \textit{``we need have no further hesitation in accepting the reason given by Vitruvius"}\cite[p.~34]{Penrose1888}. By the early twentieth century, the conflation of direct observations of curvature with Vitruvius’ words had led to its widespread acceptance as established fact. For instance, in the 1905 Frederick Moore Simpson's \textit{History of Architectural Development}, we find \cite[p. 92]{Simpson1905} ``\textit{The lines of the entablature are often not straight but rise toward the centre in a convex curve; because long lines, when quite straight, appear to `sag' or drop in the middle.}"

There are two important cases to distinguish here depending on the presence or absence of \textit{modifiers} to the main horizontal lines, such as other line segments crossing the horizontal or shapes surrounding the element. Without distractions, the situation is quite simple and well understood in vision science: in the absence of modifiers, a straight horizontal line appears straight and horizontal, and a  curved one appears  curved as  discussed below.

\subsection{The perception of curvature}

\begin{figure}[ht!]
 \centering
\includegraphics[width=0.8\linewidth]{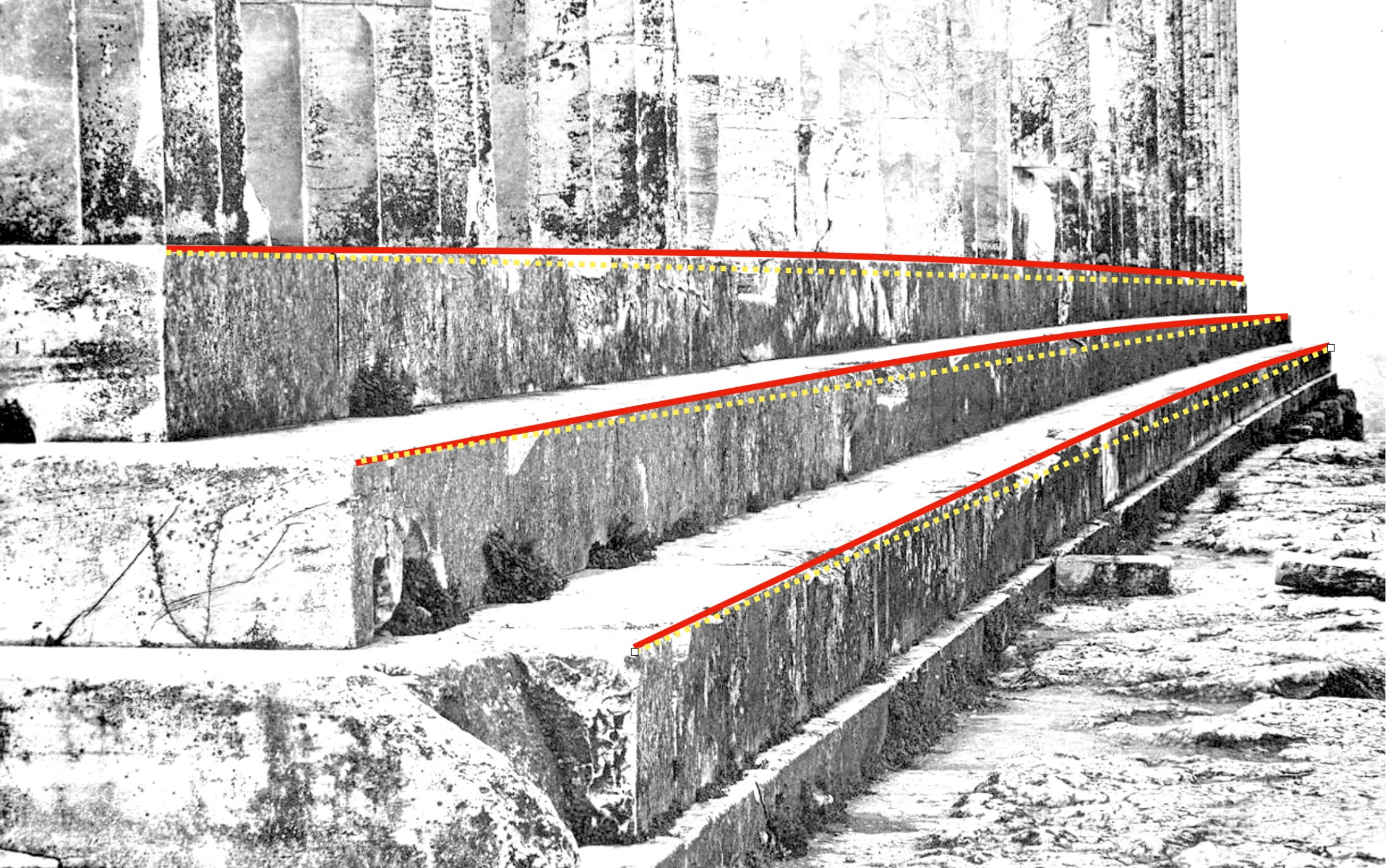}
 \caption{Viewed from the north-west corner, the stylobate appears clearly curved (1910 photograph by Fr\'ed\'eric Boissonnas \cite{Fougeres1910}). \cor{While photographs may introduced distortions, the appearance of a strong curvature viewed from the corners as been widely reported.} The dotted yellow curves are straight lines between endpoints and the red curves are the best visual fit. }
 \label{Fig-NW-corner}
\end{figure}

A simple question is  whether a slight upward curvature is  perceptible, \textit{in principle}. That is whether or not our visual system would be able to detect it at all. If not, there is really no point in discussing it. People who have visited the Parthenon know that the curvature of the stylobate can be seen from  the corners as shown in \cref{Fig-NW-corner}. Tourists are often taken to these corners so that they can be told about of optical corrections. Hence, there are  vantage points from which the curvature can be seen, which is already a contradiction regarding the possible reason behind the adjustment. What if the observer faces the Parthenon so that the two corners of the East facade are at equal distance and they see the entire span of the stylobate as in \cref{Fig-east}, say at a distance of 25 metres?

\begin{figure}[ht!]
 \centering
\includegraphics[width=\linewidth]{Fig5.pdf}
 \caption{The east side of the Parthenon from a photograph taken between 1872 and 1875. The red lines are straight and horizontal showing that at this resolution and from this viewpoint the small curvature of a horizontal element in real conditions may be difficult to perceive, \cor{with the caveat that photographs  may also introduce distortions}.  Notice also that from this viewpoint it is difficult to evaluate the thickness of the corner columns (see Myth \#3). }
 \label{Fig-east}
\end{figure}

It turns out that in ideal conditions, the human visual system is remarkably good at detecting curvature. Deviations from a straight line to a segment of a circle with the same endpoints is measured by the \textit{sagitta} $s$ in geometry, defined as the maximum perpendicular distance between a straight chord and the arc of a circle spanning the same endpoints. In experimental psychology, it is usually reported in angular units $\theta$, because visual discrimination depends on retinal angle not absolute length, and is expressed in radians (or arcmin/arcsec  with one arcmin=${\pi}/({180 \times 60} )\ \text{radians}$).  At a distance $D$, the two are simply related by $s=D\tan(\theta)$. 

Depending on the conditions and the stimuli,  studies have reported that we  can detect the largest departure from a perfectly straight line when the sagitta is as small as 3-18 arcsec of visual angle on the retina,  which, remarkably, is below the size of foveal photoreceptors
\cite{watt_ward_casco_1987_deviation_from_straightness}. However, more recent studies under different conditions give a larger estimate of 120-420 arcsec \cite{kramer_fahle_1996_detecting_low_curvatures}. If $\theta_\text{crit}$ is the threshold in radians, we can  detect at a distance $D$ a curvature due to a circle if the sagitta has height  larger than $s_\text{crit}=D\tan(\theta_\text{crit})$. Therefore, with this more conservative estimate, the human visual system could detect a line bending away from straightness by 1.5-5 cm from a 25-meter viewing distance. Since the maximal changes in elevation of the stylobate is about 6 cm and its profile well approximated by a circle or parabola \cite{Georgopoulos2012}, we could in principle detect this curvature if the observer is close enough. However, further away, it soon becomes imperceptible.

\begin{figure}[ht!]
 \centering
\includegraphics[width=0.8\linewidth]{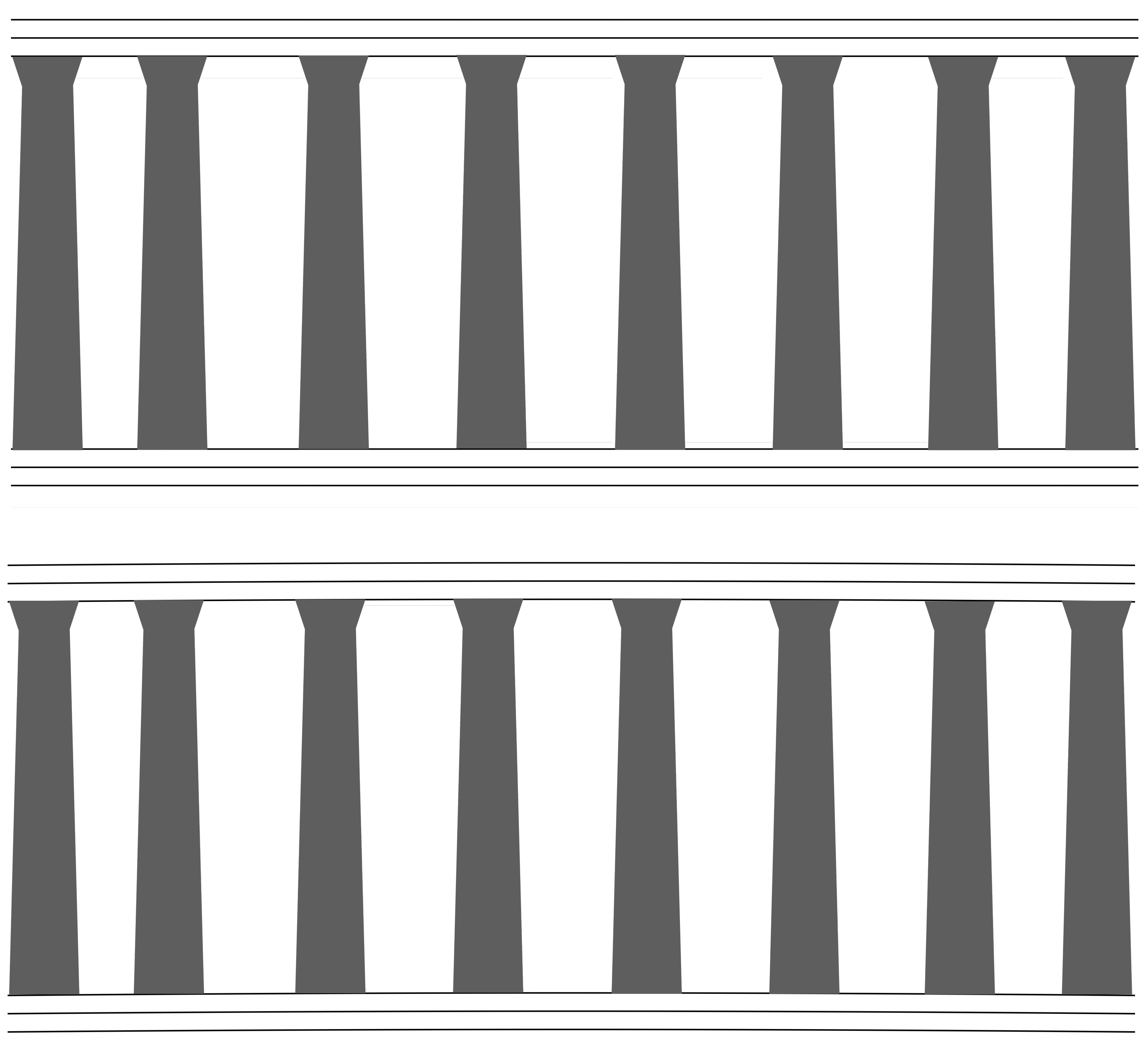}
 \caption{The comparison between two profiles: one has perfectly straight horizontal elements and one  has the exact curvature  found in the east facade of the Parthenon. It demonstrates that the horizontal lines do not appear as sagging and that an upward curvature can be detected in principle under ideal conditions. Note that the inclination of the column is so minute that on an A4 piece of paper, the virtual line extending the corner columns would meet at about 10 metres above your head. The exercise of determining which one is curved is left to the reader. }
 \label{East-ideal}
\end{figure}

To test this prediction, \cref{East-ideal} presents two two-dimensional idealisations of the east facade of the Parthenon with outlines of columns based on \cref{Fig-east}. One of this figure has  perfectly straight horizontal elements, whereas the other has the measured curvature of the stylobate taken from \cite{Georgopoulos2012}. It is not obvious from the figure to spot the curved one but it is certainly within perception limits. Not only is the human visual system very good at seeing curvature, it is remarkably good at detecting straight horizontal and vertical lines when there are no distracting elements (such as the ones found in the Hering illusion, more about that below) \cite{howe2005perceiving}.

\cor{In ideal conditions and close enough, the stylobate curvature could be detected. However, it is important to emphasize the fundamental difference between a two-dimensional drawing such as \cref{East-ideal}, a black and white photograph such as \cref{Fig-east} (which itself can be subject to many distortions), and the actual three-dimensional viewing experience in situ.} 
There is actually no line in the photograph or the real scene, the visual system creates these lines for our convenience by edge detection between zones of different intensities \cite{livingstone2022vision} and best interpolation when a line is interrupted. There are also many distracting elements in any real-life scene. Therefore, there is  a remarkable amount of visual processing that takes place when we see a building with multiple features. In particular, the interpolation for missing information is actually done by the visual system by a process well approximated by minimizing the total curvature as proposed by Mumford \cite{mumford1992elastica,petry2012perception} which would naturally flatten small irregularities.  Finally, our experience tells us that most \cor{classical or neo-classical} buildings are rectilinear. Hence, there is an expectation of linearity when looking at \cor{such} buildings. The net result is that even if upward curvature could be detected in principle, it is not easily detected in actual conditions unless it is sufficiently large and obvious. \cor{It would be informative to design a controlled experiment in which observers are asked to report whether specific architectural elements appear straight or curved when viewed from varying distances. Such a protocol would allow a direct assessment of perceptual thresholds in conditions relevant to the Parthenon. If geometrically straight elements are consistently perceived as straight across these viewing conditions, this would provide additional empirical support against the necessity of introducing optical corrections.
Since, so far, no such perceived curvatures have been reported in straight elements, we can reasonably assume that perfectly horizontal elements do not appear to sag and slightly curved ones, if perceived at all, would appear to bow upward not downward, further contradicting the original claim.}

\subsection{The illusion of curvature}

\begin{figure}[ht!]
 \centering
\includegraphics[width=\linewidth]{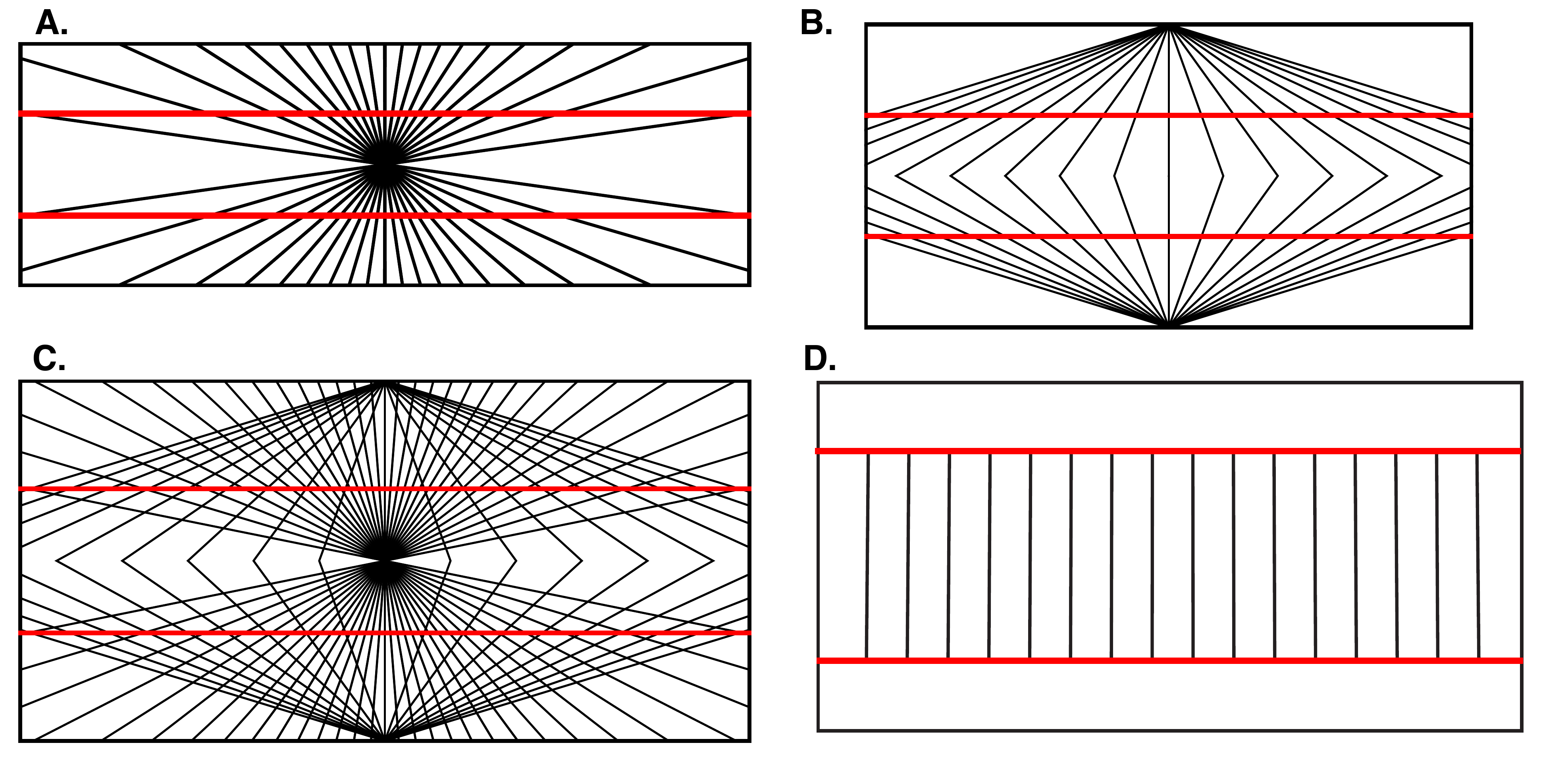}
 \caption{In the classical Hering (A) or Wundt (B) illusions straight horizontal lines appear bowed. The inversion of the effect and the nullification of the effect when both modifiers are present (C) indicate that the illusion depends on the angles of intersection   between the red lines and the modifiers. In particular when modifiers  intersect almost perpendicularly the red line (D), the illusion disappears (in this last figure, the lines have the same inclination as the centreline of the Parthenon columns).}
 \label{Fig-wundt-hering}
\end{figure}

While the visual system is very good at perceiving a straight horizontal line in the absence of distractions, the situation is different when  \textit{modifiers} are added to an image.  The most famous of such illusions is the \textit{Hering illusion} shown in \cref{Fig-wundt-hering} and first described in 1861 by the German physiologist Ewald Hering \cite{hering1861ortssinn}. When two horizontal straight lines are modified by small lines, they appear to bow. The inversion of the modifiers lead to the opposite effect, the \textit{Wundt illusion}\cite{wundt1898geometrisch}, and the two effects together nullify the illusion. Hence, we  conclude that it is the angle between modifiers and the horizontal line that is relevant for the illusion. This type of illusion is consistent with \textit{Brentano’s Law},  the visual tendency to overestimate  acute angles (as larger than they actually are), and underestimate  obtuse angles (as smaller than they are in reality), a phenomenon that is well confirmed experimentally \cite{howe2005perceiving}.

\begin{figure}[h]\begin{center}
\includegraphics[width=\linewidth]{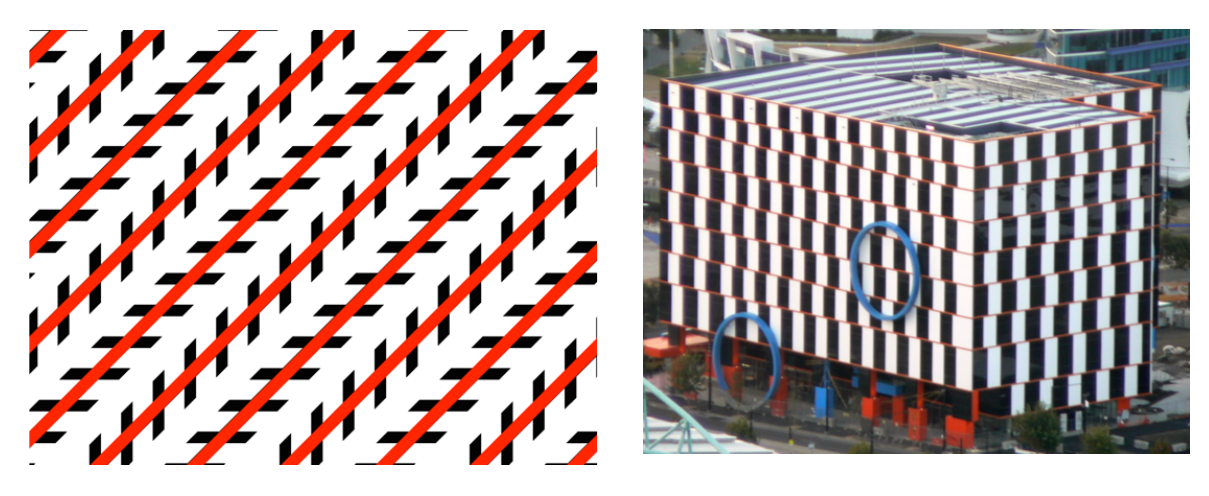}
\caption{In the Z\"ollner illusion, parallel lines (in red) appear to be converging or diverging. A variation of the same illusion is used on  the Port 1010 building (Docklands, Melbourne, Australia. Photograph courtesy of Joe Bekker) }
\label{gr6-zollner}\end{center}
\end{figure}
Another classic illusion is  the\textit{ Z\"ollner illusion}, discovered in  1860 by Johann Karl Friedrich Z\"ollner, a German astrophysicist, who noticed that parallel lines appeared to diverge or converge when crossed by short diagonal lines \cite{zollner1860pseudoskopie}. Though the main lines are perfectly parallel, they look skewed or tilted because of the added intersecting strokes as shown in \cref{gr6-zollner}. Variations of this illusion such as the Caf\'e-wall Illusion have been used in actual constructions (see \cref{gr6-zollner}). But note that the curvature is not affected in this particular illusion.

Before we proceed with the discussion of illusions, it is important to note that the perception of visual illusions is very subjective with large variations in amplitude between individuals depending on many individual characteristics (e.g. age, gender, eye corrections, prior exposure to illusions,...) and many external factors (e.g. light, color, distance, orientation of the stimuli, orientation of the head,...).\footnote{This variability implies that you, the reader, may not perceive the visual illusions described in this paper, which is no reason to worry. It just means that you have a good visual system.} The important implication of this variability for the discussion is that there will never be  one-size-fits-all correction for any visual illusion. It is easy to change the curvature of the red lines in \cref{Fig-wundt-hering}A so that they appear again straight, but this correction will depend on many different factors even for a single individual.

A natural question is whether such curvature illusions could be at play in buildings and if it was somehow used in the Parthenon. The first thing to realize is that the illusion completely disappears if the modifiers intersect the object at perpendicular angles as shown in \cref{Fig-wundt-hering}D. Even close to the perpendicular, no effect of curvature is induced. This is attributed to the fact that the visual system is much better at detecting horizontal, vertical, and perpendicular lines.  Hence, the illusion disappears as the perceived angle is the very close to the actual angle and the visual system does not have to curve the object to resolve a possible disagreement. Since the main modifiers are the columns and they are straight, we conclude that this illusion does not take place.

A first possible rebuttal to this last argument is that the columns have a slight inward tilt so that their centrelines all converge to a point located at a middle point in the plane of the columns above the stylobate. The estimation of this angle is around  0.4$^\circ$ \cite[Fig.~3.29]{Korres1998}.
which means that the point of convergence for the east facade is about  2.2 km above the stylobate. This angle is so minute that, even in ideal conditions on paper, it does not induce curvature as can be appreciated from the absence of such an illusion in \cref{Fig-wundt-hering}D. Furthermore, these inclinations are much smaller than the angles formed by the edges of the tapered columns, which dominate the visual scene. However, since the edges of the tapered columns have two sides with opposite inclinations, their effect on a possible Hering illusion cancel each other, very much like in the Hering-Wundt illusion of \cref{Fig-wundt-hering}C. 

A second rebuttal is that in a three-dimensional environment, the lines would not appear perpendicular when viewed from below, since their projections shift toward a vertical \textit{vanishing point}. This explanation would at least generate an illusion in the expected direction, as illustrated by the lower line in \cref{Fig-wundt-hering}A. However, several issues arise with this line of reasoning. First, the inclination depends on the viewing point. It is strong at very short distances and gradually disappears with the distance. It would therefore be impossible to correct such a possible effect for all viewing distances. Second, humans have a compensatory mechanism for perspective. When we look at the Petronas Twin Towers in Kuala Lumpur from the ground they do not appear vertical but inclined towards each other.
However, nobody interprets this image as being two towers leaning towards each other. We naturally correct for this perspective and the skybridge between the towers is not perceived as curved because the two vertical lines are slightly inclined on their retinal projection. 

According to the theory of Howe and Purves \cite{howe2005perceiving}, two-dimensional paper illusions such as Hering's occur precisely because our visual system is best trained to evaluate angles in the real three-dimensional world. The errors in evaluation that we make on paper are due to the fact that our visual system is tuned for real-world vision. Hence, there is no reason to expect them in viewing the Parthenon.

\begin{figure}[h]\begin{center}
\includegraphics[width=0.8\linewidth]{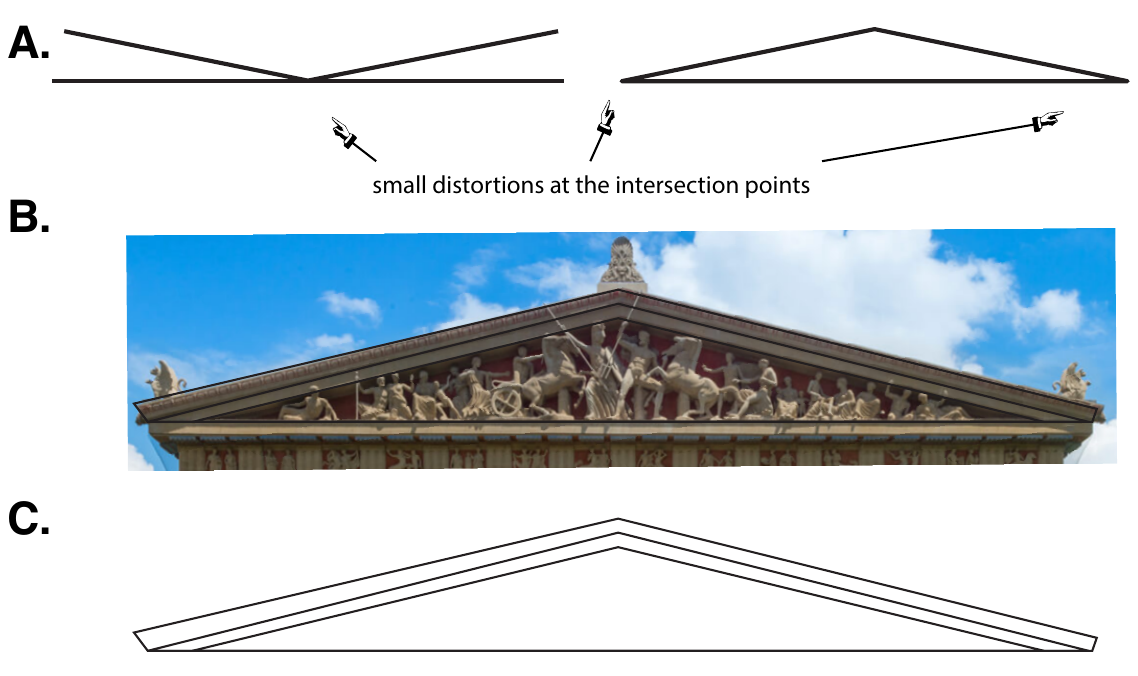}
\caption{A. The intersection between a horizontal line and another line at an acute angle creates small distortions to the horizontal line close to the intersection. B. The triangular pediment of the Parthenon's replica in Nashville, Tennessee (originally built in 1897 with the same architectural refinements). C. The outline of the  pediment does not show sagging in the horizontal element.}
\label{Fig-pediment}\end{center}
\end{figure}

A last possible rebuttal in favor of a visually-induced curvature is found not in the stylobate but in the triangular pediment. This is the only place in a Greek building with obvious elements that are neither primarily horizontal nor vertical (disregarding the small column inclination, their tapering and curvature --see below). There the argument goes that the oblique edges of the triangular structure could induce a Hering-like illusion. This argument first put forward by August Thiersch in 1873 \cite{Thiersch1873} rests on the same type of illusions where small modifiers on a straight line change their orientation as shown in \cref{Fig-pediment}A. However, the small distortions created by the same Brentano effect are localized at the intersection points as first described by Hermann von Helmholtz in his treatise on optics \cite{von1867handbuch,von1925helmholtz} and discussed in \cite[p.~218]{sanford1903experimental}. If we extract a triangular outline from the pediment (\cref{Fig-pediment}B), no downward curvature is observed not even in its ideal two-dimensional shape (\cref{Fig-pediment}C). This illusion, used to justify the curvature of the architrave, does not exist. 

\subsection{Alternative theory: Curvature for rain drainage}
An interesting suggestion for the curvature of the stylobate is that it was conceived for drainage.
The idea is simple. If the stylobate is domed, it would naturally help rainwater drain away from the platform.
The Parthenon was roofed with a timber structure carrying overlapping Pentelic marble tiles with channels to guide water to the eaves. This roof spanned only the inner cella and its porches, protecting the cult statue and interior walls from direct rainfall and leaving the surrounding peristyle  open: the colonnaded area standing on the stylobate was not roofed. Thus large portions of the stylobate were  exposed to rain. The beautiful lion-heads at the roof corners were decorative pseudo-spouts, not functional drains. Therefore, rainfall striking the stylobate ran outward across the marble steps. 

We can estimate, based on the stylobate profile, how effective the curved platform would have been to prevent pooling near the temple base. From the profile (\cref{fxy}), we  average the slope for the four sides and obtain an angle of $\theta_\text{ew}\approx 0.24^\circ$ on the short east and west side and $\theta_\text{ns}\approx 0.18^\circ$  on the long north and south sides (note that since the central cella was roofed, only the peristyle is relevant here).

The Parthenon, including its stylobate, was constructed from {Pentelic marble}, a fine-grained calcitic marble that is nearly impermeable, with very low porosity. Hence, we can assume that water does not infiltrate into the platform but runs over its surface.

A thin, laminar film of thickness $h$ flowing down a plane at slope $\theta\ll 1$ is governed by  the classical Nusselt film solution \cite{nusselt1916,ruyerquil2000,bird2002} that gives a parabolic velocity profile  with depth-averaged velocity $U$ and flux per unit width $q$ given by :
\begin{align}
U &= \frac{g\sin\theta}{3\nu}\,h^{2}, \qquad 
q = Uh = \frac{g\sin\theta}{3\nu}\,h^{3},
\end{align}
where $h$ is the film thickness, $g$ gravity, and $\nu$ is the kinematic viscosity of water. It is reasonable to assume that the well-worn marble enforces a no-slip boundary condition, so these relations hold without additional correction (which would have the effect of pinning the film and slowing down drainage). The formation of rivulets is also possible but does not change significantly the estimate.

At $20^\circ\mathrm{C}$, we have for water $\nu\simeq 10^{-6}\,\mathrm{m^{2}\,s^{-1}}$. With local slopes $\theta\approx 0.2^\circ$ and a film of $h=1$ mm, we obtain  $ U\approx 0.3\, \mathrm{m\,s^{-1}}.$  Over a maximal  runoff distance of $L=30\,\mathrm{m}$, a characteristic time is $t=L/U\approx 35$ minutes. But since the film decreases as it drains, the drainage time for an initial film of thickness $h_0$ is 
\begin{align} t_{\text d}=\frac{L}{\langle U \rangle},\quad  \text{where}\quad \langle U \rangle = \frac{g \sin \theta}{3\nu} \, \langle h^2 \rangle=\frac{g \sin \theta}{9\nu}  h_0^2,\end{align}
which gives $t_{\text d}=2.2$ hours for $h_0=1$ mm. It would then take  hours for the water to drain over the full span. Therefore, we conclude that the small gradients help water to drain slowly and prevent the formation of puddles on the stylobate. Whether it was designed with that purpose in mind is unknown, but we find that Penrose's dismissal of this effect unjustified. However, this does not explain why the curvature also appears in the architrave.

\section{Myth \#2: Column entasis}

\subsection{The refinement}
\textit{Entasis} is the slight swelling (concave) of a column shaft that departs from a purely straight linear tapering. It appears in the early sixth century BCE in Doric temples, such as the first temple of Hera at Paestum (ca. 550 BCE) where   this curvature is very pronounced and obvious to the eye. By contrast, in the Parthenon, the entasis is so subtle that it escaped notice for centuries and was only rediscovered in 1814 by Charles Robert Cockerell during his architectural investigations in Athens \cite[p.~48]{pearce2017charles}\cite{salmon2008cr}. The word \textit{entasis} was first  used by Vitruvius and derives from the Greek word $\mathrm{\acute{\epsilon}\nu\tau\epsilon\acute{\iota}\nu\omega}$ (enteino), associated with "tension" or "to stretch". Although mostly attributed to Greek architecture, evidence indicates that the principle of entasis was already present in Mesopotamian and Egyptian architecture centuries earlier. Unlike the curvature of the stylobate, the small columnar curvature did not disappear after the Doric period. Indeed, Roman builders adopted entasis so enthusiastically and universally that no Roman column without entasis has ever been reported \cite{stevens1924entasis} and the convention was carried forward without interruption into the Middle Ages. With the Renaissance revival of Vitruvius’ treatise, interest in entasis was renewed in Italy and rapidly disseminated throughout Europe, eventually reaching the Americas through neoclassical architecture.

Since Vitruvius' method and drawings for creating a column with entasis were lost, a long-standing problem, first formulated in Nicolas-François Blondel's ``\textit{Résolution des quatre principaux problèmes d’architecture''} (1673), is the exact mathematical nature of the profile and how it was designed \cite{raynaud2020mathematiques}. The list of possible mathematical functions is long and includes segments of the parabola, the circle, the ellipse, the hyperbola, the catenary, the conchoid of Nicomedes, and the vertical projection of a helix \cite{stevens1924entasis,pakkanen1997entasis,radelet2025entasis}. The discovery of  third-century sketch drawings incised into the unfinished walls of the Temple of Apollo at Didyma suggests that an initial circle stretched to an ellipse may have been the method used \cite{haselberger1985construction}.

\begin{figure}[h]\begin{center}
\includegraphics[width=\linewidth]{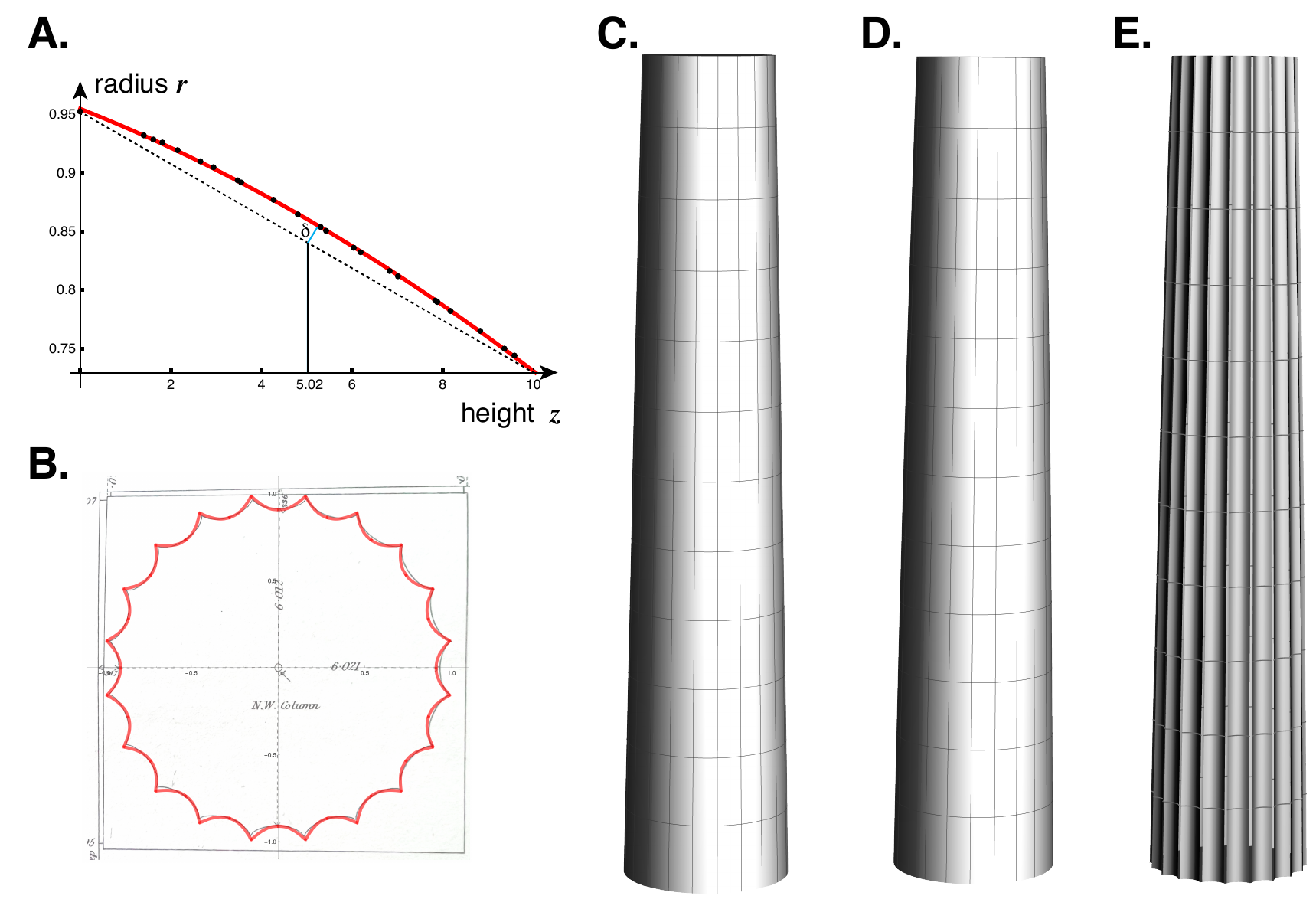}
\caption{A. The radius as a function of height is well approximated by a parabola (based on Penrose's measurements \cite[Plate 14]{Penrose1888} of the north-side central column of the east facade). Here $r(z)=0.954753 - 0.0152412 z - 0.000716763 z^2$. A linear tapering would have been represented by the dashed line. The maximal deviation from a linear tapering is attained close to the middle of the columns (at 5.02m with a deviation of $\delta\approx18$mm). B. The fluting of Doric columns is well approximated by a simple power law ($\alpha=6.4$ and $\beta=2.25$, cross-section of the north-west column from \cite[Plate 13]{Penrose1888}). CD. The reconstructed column with and without the measured entasis. The exercise of determining which one is curved is left to the reader. E. The same column with flutes.}
\label{Fig-entasis}\end{center}
\end{figure}

For our study, since the swelling is so minute, the true function and its design are not directly relevant. Indeed,  the profile of the Parthenon columns is  well captured by a simple parabola with three parameters (whereas a linear tapering requires two parameters) as demonstrated for a single column in \cref{Fig-entasis}A. Denoting $r$ the radius of the column from its centerline (assumed to have an axial symmetry) and $z$ the height, we have simply
\begin{align}
    r(z)=c_0+c_1 z+c_2 z^2,\qquad  z\in[0,h],
\end{align}
where $h$ is the column height.

Another obvious characteristic of the columns is the flutes around the columns. These vertical grooves were cut into the shaft of a column after it was erected and run seamlessly along the entire column, representing a major construction challenge for the 46 outer and 23 inner columns of the Parthenon. Each column  carries exactly 20 flutes  with sharp arrises, a  characteristic feature of the Doric order. 
For modeling purposes, the unit function describing a single flute can be chosen to be a power law (starting at the middle point of a groove):
\begin{align}
    f(\theta)=\frac{1+\alpha\, |\theta|^\beta}{1+\alpha\, (\pi/10)^\beta},\qquad  \theta\in[-\frac{\pi}{20},\frac{\pi}{20}],
\end{align}
where the scaling implies that the maximum of $f$ is equal to  1 at $\theta=\pm\pi/20$.
Then, a model for the entire column is given by a  function with five parameters describing the radius at height $z$: 
\begin{align}
    F(\theta,z)=r(z) \text{Per}(f(\theta),\pi/10),\qquad  \theta\in[-\pi,\pi],\ z\in[0,h],
\end{align}
 where Per$(f(\theta),\pi/10)$ is the $\pi/10$-periodic extension of the function $f(\theta)$. Examples with and without refinement, and with and without flutes are given in \cref{Fig-entasis}.

\subsection{The myth}
\begin{figure}[ht!]
 \centering
\includegraphics[width=0.8\linewidth]{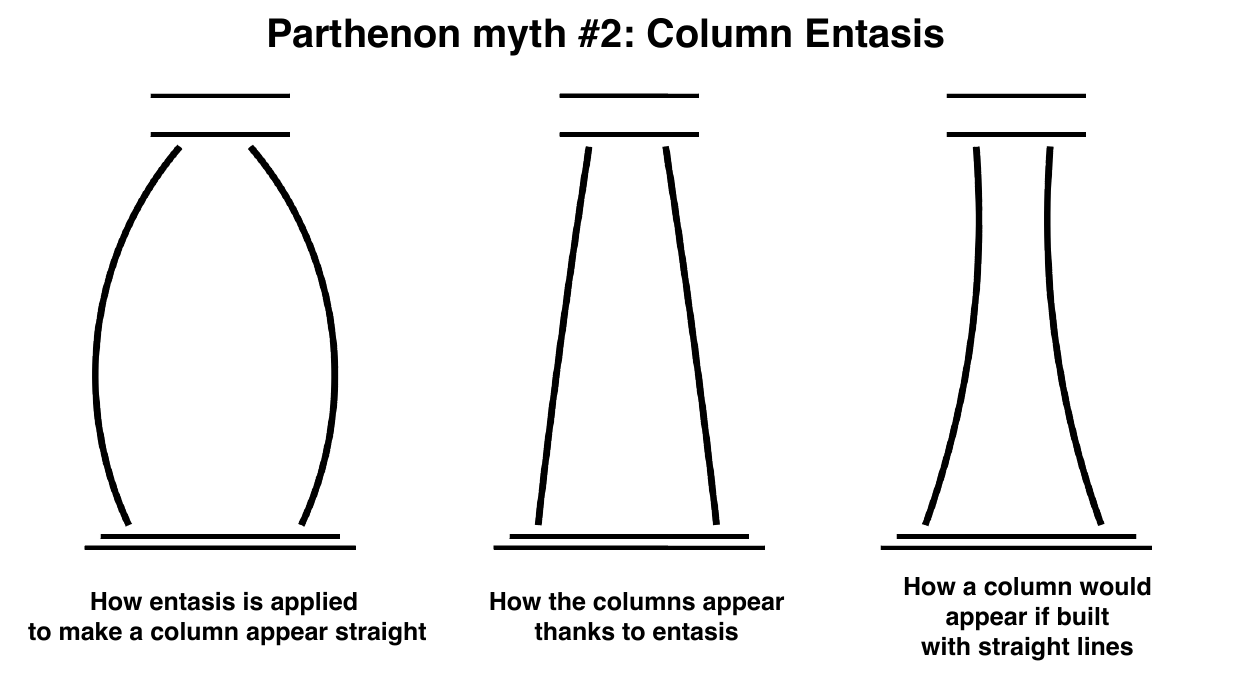}
 \caption{The myth of the entasis explains the small swelling of columns as a correction to an optical illusion. None of these statements is correct.}
 \label{fig2}
\end{figure}

The myth related to the entasis is clearly expressed by Penrose \cite{Penrose1888}:
``\textit{The entasis  is the well-known increment or swelling given to a column in the middle parts of the shaft for the purpose of correcting a disagreeable optical illusion, which is found to give an attenuated appearance to columns formed with straight sides, and to cause their outlines to seem concave instead of straight. The fact is almost universally recognized by attentive observers}''. In simple terms, ancient architects added a small curvature to the column profile so that they would look straight as shown in \cref{fig2}.

The myth is usually justified by two sentences found in Section 3.3.13 of Vitruvius.\footnote{``\textit{Haec autem propter altitudinis intervallum scandentis oculi species adiciuntur crassitudinibus temperaturae. Venustates enim persequitur visus, cuius si non blandimur
voluptati proportione et modulorum adiectionibus,
uti quod fallitur temperatione adaugeatur, vastus et invenustus conspicientibus remittetur aspectus.
''}}
It is important to recognize that translating Vitruvius is by no means straightforward, and some translators have approached his writings with a preconceived notion of optical corrections. A recent translation reads \cite{VitruviusRowland1999}: ``\textit{These adjustments to the diameter are added because of the extent of the distance for the ascending
glance of our eyes. For our vision always pursues beauty, and if we do not humor its pleasure by the proportioning of such additions to the modules in order to compensate for what the eye has missed, then a building presents the viewer with an ungainly, graceless appearance."} While clearly referring to refinements to please the eye, there is no direct reference in this translation of Vitruvius's writings to the illusion of a column with straight edges appearing hollowed, despite this interpretation found in other translations (e.g. in Granger \cite[p.~180]{vitruvius1931} ``\textit{Without a slight swelling of the shaft of a column, the straight upright line would strike the eye as hollowed inwards.''}).

\subsection{The illusion of an illusion}

This myth is particularly perplexing. First, its foundation is wholly unfounded: the assertion that straight-sided shafts appear concave lacks empirical support, and no such illusion has ever been scientifically demonstrated, even in controlled tests using three-dimensional reconstructions of temples directly aimed towards this question \cite{thompson2007origins}. Second, in the Parthenon the entasis is so slight that, according to the estimate from the previous section and an 18 mm sagitta, the curvature is below our ability to be perceived from most distances (a test is proposed to the reader in \cref{Fig-entasis}). Moreover, when perceived, it  registers to the eye as a bulging  rather than as a correction. Third, while entasis is pronounced in many other temples before and after Pericles's period, it is extremely subtle in the Parthenon.   Are we to assume, then, that architects initially introduced entasis as a corrective measure, but executed it ineffectively for centuries, only to discover much later how to apply it properly? That they only applied it correctly for the Parthenon and  in a way that  may not even be perceived. And that later on, entasis would be used by Romans for centuries in a very prominent way that indicates they did not understand the reason behind this adjustment? This line of reasoning strains credibility.

A simpler alternative explanation is that builders used curved columns because they believed them more pleasing than straight ones. 

\subsection{An alternative theory: The organic explanation}

Another possible explanation for refinements is that every refinement is seen as an expression of a living, quasi-biological vitality, a manifestation of a generative force analogous to growth in nature,  a tree that swells upward from its roots, or to a flexed muscle in tension all governed by curves and rhythm.  \cor{It has been argued that this type of vitality  appears in Gothic architecture as described by Ruskin in the \textit{Stone of Venice} \cite{ruskin1867stones}. In the context of Greek architecture, these curves are taken to endow the temple} with a sense of life, elasticity, and organic embodiment in contrast to the ``vulgar" rigidity of purely rectilinear constructions \cite[p.324]{Choisy1899}. There is great poetry associated with this “organic theory” and many authors have unleashed their creative writing skills along these lines. 
The great advantage of such an organic interpretation, apart from resonating with modern sensibilities, is that it is, as Karl Popper would say, completely immunized against criticism. Indeed, this type of explanation is metaphysical and  claims about vitality, organic energy, or natural flows cannot be tested, falsified, or verified. It sits in the mind of the beholder and, from a scientific point of view, nothing more can be said about it.

\subsection{Another alternative theory: column buckling}

Yet another explanation for entasis is that it improves the resistance of columns against buckling. Interestingly, this is yet another enduring myth that originated soon after the initial development of a theory of buckling by Euler in 1744 \cite{eu44,eu59}. In his 1770 treatise \textit{``Sur la figure des colonnes''} Joseph-Louis Lagrange directly addressed the problem of entasis \cite{la70}. He criticized the fact that the shape used has no proper justification and attempted to find, through mathematical analysis, the shape of a column of fixed volume and height that can bear the largest compressive load.  He concluded, erroneously, that the optimal shape is a straight cylinder.

It is  well known nowadays that an \textit{elastic} column under load can become unstable through Euler buckling  \cite{goriely2008nonlinear,dedego11} and that the critical load can be increased  by modulating the cross section \cite{keller1960shape}. 
Therefore, it is no surprise that some authors have tried to justify entasis as a means to increase the stability of marble columns \cite{cox1992shape}.  While this type of analysis is particularly important when applied to trees or steel beams,  it is not directly applicable to Greek columns \cor{or any type of masonry  for that matter \cite{heyman1966stone}}. Indeed, Greek columns were typically constructed from stacked marble drums, connected by metal dowels or wooden pegs.  Their failure mechanisms have been the subject of structural and seismological studies, since ancient ruins often exhibit collapsed peristyles \cite{PapaloizouKomodromos2011,DrososAnastasopoulos2015,VassiliouMakris2012}. These studies found that archaeological evidence and experimental shake-table studies consistently show that multi-drum marble columns fail through the following mechanisms: (i) \textit{Brittle fracture of marble} (i.e. cracking due to the low tensile strength of marble); (ii) \textit{Rocking and sliding at drum joints} (main failure mode under seismic actions); (iii) \textit{Overturning} (i.e. progressive collapse during earthquakes); and (iv)  \textit{Foundation settlement} (leading to eccentric loading and cracking). Elastic buckling is notably absent from observed collapse modes. 

The fact that buckling is not relevant for Doric columns is further demonstrated by a simple estimate. The Euler critical load for a prismatic column is \cite{goriely17}
\begin{equation}
    P_{\mathrm{cr}} = \frac{\pi^3 r^4 E}{4 L^2},
\end{equation}
where $E$ is Young's modulus, $L$ is column height, and  $r$ is its radius. For  typical Pentelic marble $E$ is around 60 GPa, but a study on the degradation of Pentelic marble from the Acropolis monuments shows a reduced value\cite{Skoulikidis1993}. Hence, we take conservatively $E \approx 40$ GPa, $L \approx 10.4$ m, $r \approx 0.75$ m (corresponding to the smaller radius at the top of the tapered columns), to obtain a lower bound for the Euler critical load and average  buckling stress:
\begin{equation}
    P_{\mathrm{cr}} \approx 
    0.907 \times 10^9 \ \text{N},\qquad 
    \sigma_{\mathrm{cr}} = \frac{P_{\mathrm{cr}}}{\pi r^2} \approx 513 \ \text{MPa}.
\end{equation}
 This last value exceeds the compressive strength of Pentelic marble  which is around 100 MPa \cite{Ruedrich2013}, with other estimates for marble as low as 50 MPa.  The fact that  marble columns would crumble before they buckle is mainly due to the fact that Doric columns are not slender, which is a key ingredient for the occurrence of elastic buckling.
 We conclude that  crushing or joint sliding occur long before Euler buckling, which is therefore irrelevant to the discussion of its shape.

\section{Myth \#3: Larger corner columns}
\subsection{The refinement}

Another example of irregularities found in the Parthenon is that the corner columns are  thicker (around 1.947 m at the base) than the flank and central columns (with an average closer to 1.905 m) \cite[p. 13]{Penrose1888}. This 4.2 cm difference is on the order of 2\%, which may just be visually perceptible at the scale of the colonnade. All four corner columns are made thicker, while the diameters of the others remain more uniform. Again, this difference suggests a deliberate choice.

\subsection{The myth}
The explanation for this increase in thickness as a correction to an illusion appears explicitly in Vitruvius  and all translations are consistent. For instance \cite[3.3.11]{VitruviusMorgan1914}  ``\textit{The columns at the corners should be made thicker than the others by a fiftieth of their own diameter, because they are sharply outlined by the unobstructed air round them, and seem to the beholder more slender than they are. Hence, we must counteract the ocular deception by an adjustment of proportions.''}

\begin{figure}[ht!]
 \centering
\includegraphics[width=0.8\linewidth]{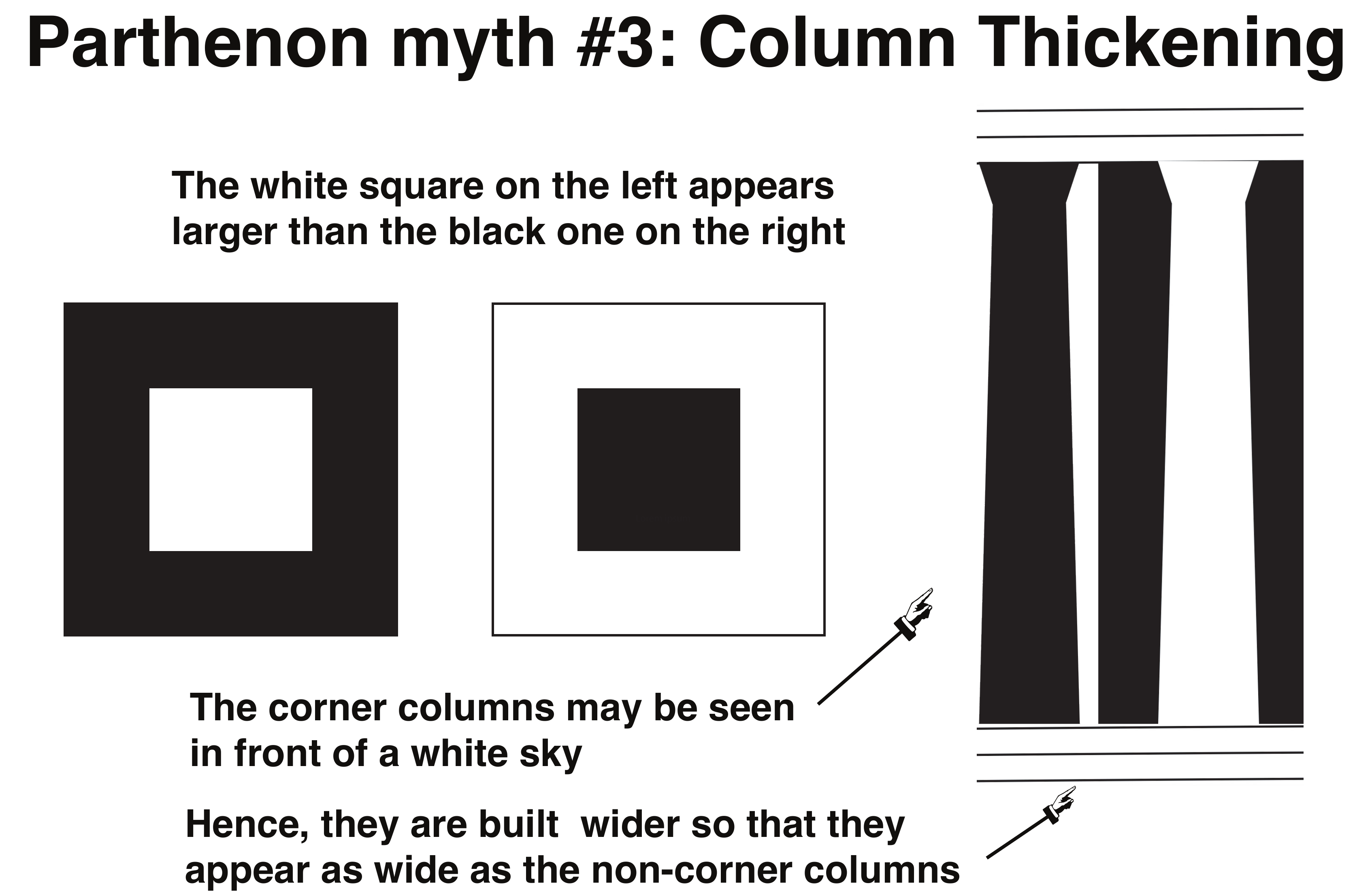}
 \caption{The irradiation illusion. The white square on a white background appears slightly larger than the black square on a white background. The myth associated with  corner columns is that they were designed thicker to compensate for this illusion.}
 \label{Fig-irradiation}
\end{figure}

In modern times, the illusion that has been identified as a possible explanation is the so-called \textit{irradiation illusion} \cite{plateau1873statique} in which a white object on a dark background appears bigger than  a dark object of the same size on a white background as seen in \cref{Fig-irradiation}.
An explanation of this illusion is that the white light from the brightly lit object stimulates adjacent regions in the visual field, blurring boundaries and causing overestimation of the bright area \cite[p.~93]{Robinson2013}. 

The extrapolation of this effect fueled the following myth: when a brightly lit column is seen against a dark background, its diameter may appear larger than it actually is (see \cref{Fig-irradiation}). In principle, if all columns were of the same diameter, the corner column would then look thinner by comparison if seen against a bright sky than the other columns seen against a dark interior. The myth is again that to correct this effect, the architects of the Parthenon somehow may have been aware of this illusion and designed slightly thicker corner columns to counteract it. This is a beautiful and simple explanation that specialists in visual illusions love to refer to in order to emphasize the importance of such effects \cite{coren2020seeing}. Apart from the fact that there is no historical evidence that the ancient Greeks were aware of this illusion and that they decided to modify the corner columns accordingly, there are many problems with this myth.

The first problem is that the irradiation illusion has only been studied in ideal conditions in experimental psychology where subjects are asked to compare rectangles on paper under different lighting conditions \cite{long1984task}. However, it has been argued that in non-laboratory setting these illusions may not appear at all \cite{mather1986irradiation}. To my knowledge, the irradiation illusion has never been empirically tested in everyday life, and there is no report of corner columns appearing perceptibly thinner in the numerous colonnades found worldwide where all columns are of equal dimensions. Hence,  the existence of this illusion in real-life settings has yet to be confirmed.

The second problem is that in order for the illusion to appear, the non-corner columns must be brightly lit and be against a dark background. However, in such a case, the corner columns would also be brightly lit and either be against a dark background themselves (see below) or a bright background. Hence, the corner columns are not in a suitable light environment that would be needed for comparison. 

\begin{figure}[ht!]
 \centering
\includegraphics[width=0.8\linewidth]{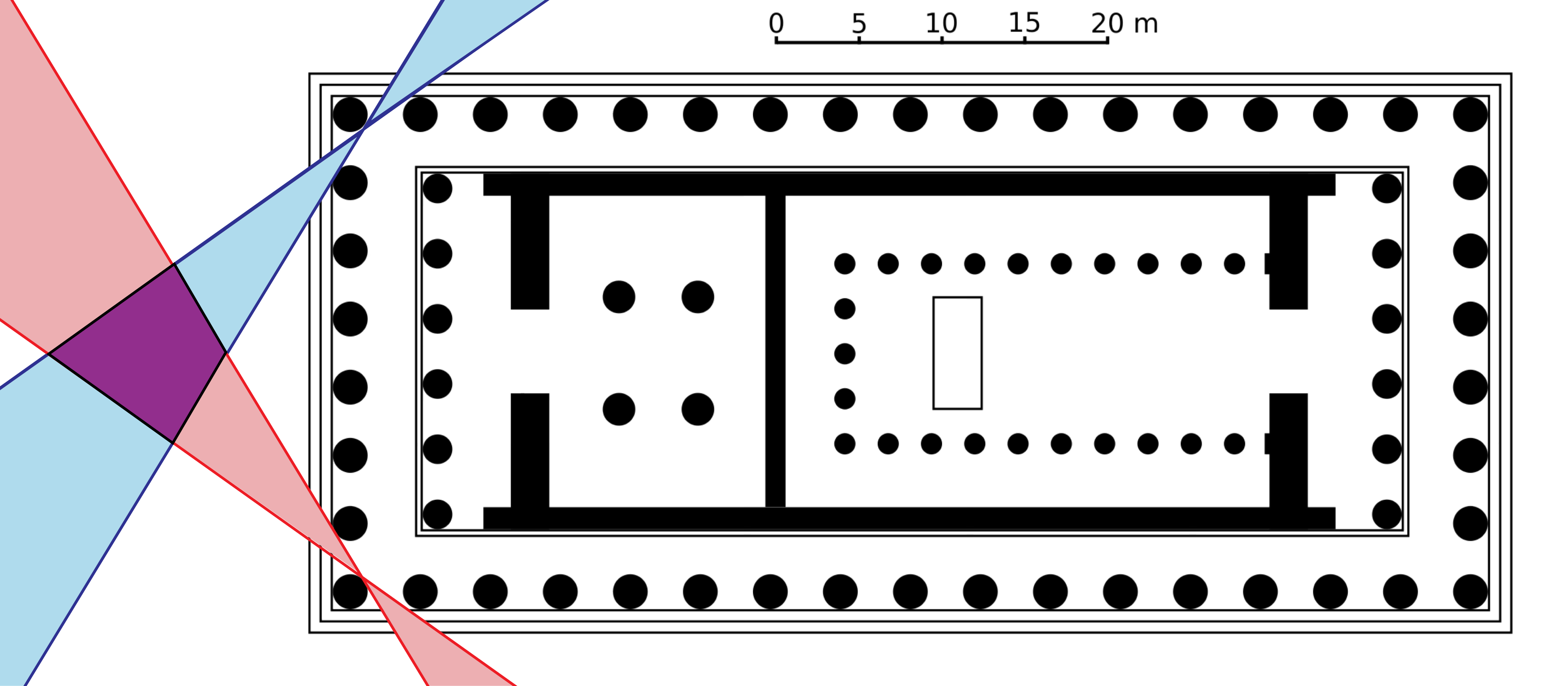}
 \caption{The blue and red areas indicate vantage points where a single corner column is visible against the sky, while the purple region marks the limited positions from which both corner columns can be  seen  simultaneously against a uniform sky. The three extra purple regions for the other facades are shown in the next figure.  The  implication of the optical-correction myth is that, from the majority of viewing positions, the corner columns would appear unequal in size relative to the others, in direct contradiction to its stated purpose.}
 \label{corner-columns}
\end{figure}

\begin{figure}[ht!]
 \centering
\includegraphics[width=0.8\linewidth]{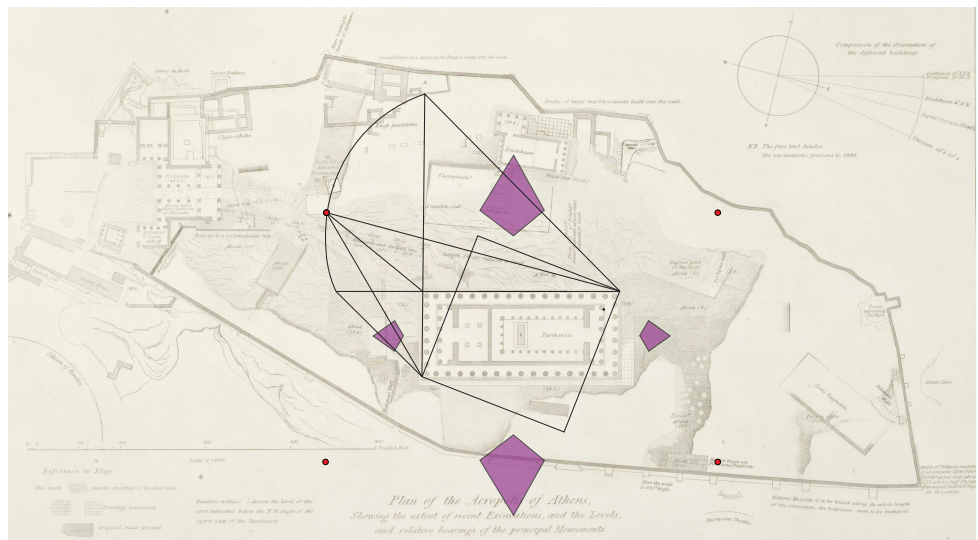}
 \caption{The purple areas on the map of the Acropolis (from \cite{Penrose1888}) indicate the special vantage points where, possibly, all the columns would be seen  of equal thickness (assuming that at these points, the lighting conditions would be perfect and that the visual system would make up exactly for the 2\% difference in column thickness, independently of the distance). The red dots correspond to the construction of Pennethorne for the points of view with the perfect harmony in proportions \cite{Pennethorne1878}. As expected, these special vantage points do not coincide with the purple areas. }
 \label{corner-columns2}
\end{figure}

The third problem is related to the viewpoint. For the explanation to make sense, the observer must be in a  place where the corner columns are against a light sky background and the other ones against a dark background. A casual inspection of the multiple  photographs of the Parthenon as the ones shown in \cref{Fig-east} show that the two corner columns are rarely seen simultaneously against a uniform white sky, with outlines typically set against different backgrounds of stone or shadow. Hence their width cannot be appreciated.  In this case, the conditions  for the irradiation illusion are not fulfilled, and any perception of the corners as visually corrected could not arise. Even if one were to assume that the illusion actually exists in real settings, its correction could only have been effective in the very restricted areas of observation shown in purple in \cref{corner-columns}, where the viewer’s angle and background alignment allowed both corner columns to be seen at once against a bright, undifferentiated field. On the map of the Acropolis shown in \cref{corner-columns2}, these four special areas  distinguish themselves by the absence of any distinguishing features  falling either outside the Acropolis, inside a building, or on a slope.
Thus, rather than producing the claimed visual equality, the adjustment would introduce irregularity across the majority of viewing conditions. If the alleged objective was to ensure the equal apparent thickness of all columns under typical circumstances and in most vantage positions, the logical and only solution would have been to employ uniform diameters throughout.

\section{Myth \#4: Proportional adjustments}

\subsection{The refinement}
When a human figure is inscribed at an elevated position, such as on a frieze, then by simple geometry, the straight-line distance from the viewer’s eye to the head of the figure is greater than the distance to the feet. The idea behind such a refinement is that sculptors deliberately alter proportions to counteract perspective effects: by enlarging the head, its retinal image becomes comparable in scale to that of a head viewed at eye level. The same is true for inscriptions above eye level. Precise measurements of the human figures on the Parthenon west frieze show a body/head ratio (BHR) ranging from 7.40 to 7.45 \cite{Bouzakis2016}. Since, standard modern anthropometric references give a corresponding ratio of about 7.95–8, the smaller BHR found in the Parthenon has been interpreted as a refinement. 

\subsection{The myth}
Could the observed BHR for the Parthenon frieze be justified from such a correction? The Canon of Polykleitos (5th century BCE) was a theoretical treatise, now lost, in which the sculptor Polykleitos of Argos set out mathematical rules for the ideal proportions (known as \textit{symmetria}) of the human body. It is not clear what this proportions were from available texts but  the best known examples of  Polykleitan statues show a BHR of 7 for the Doryphoros \cite{Tobin1975} and 7.4 for the Diadoumenos \cite{richter1935another}, close to the BHR found in the Parthenon. By contrast, the Canon of Lysippos gives a BHR of 8 according to Pliny the Elder \cite{bieber2019iii}. But Lysippos only came a century after the Parthenon and there is no evidence that his canon was in place at the time the Parthenon was built.

Further, the standard anthropometric reference with a BHR of 8 corresponds to humans in the twentieth century. Since average height has increased considerably and BHR changes with height, this reference must be adapted to ancient Greeks. Male stature at the time of the Parthenon was about 165cm compared to 180cm today \cite{koukli2023stature}. The anthropometric BHR reference for such a height is closer to 7.5  \cite{gordon20142012}. Hence there is no evidence that such a refinement is present in the Parthenon.

\section{The bigger myth}
The belief that such a scaling correction is  applied to inscriptions or body parts is widespread. The idea originates in Antiquity, first in Plato when he described  colossal works of sculptors in \textit{Sophist} (235d-236a): "\textit{If they reproduced the true proportions of their beautiful subjects, you see, the upper parts would seem smaller than they should, and the lower parts would appear larger, because we see the upper parts from farther away and the lower parts from closer.}" A similar comment is found in Vitruvius in Book 3.5.9 as translated by Raymond "\textit{To preserve a sensible proportion of parts, if in high situations or of colossal dimensions, we must modify them accordingly, so that they may appear of the size intended}" \cite[p. 253]{raymond1909proportion}. 

The Roman inscription found on the base of Trajan’s column (\textit{Corpus Inscriptionum Latinarum}, Volume VI no 960)  is often  cited as the best evidence of deliberate optical correction. Yet precise measurement shows this claim to be unfounded. Although the letter heights vary between approximately 85 mm and 115 mm from bottom to top, the variation does not follow a consistent or purposeful progression. Notably, the second line is larger than the first, contradicting the expectation of a systematic increase with height, while subsequent lines fluctuate irregularly, with lines 1–4 around 114 mm, line 5 reduced to 105 mm, and line 6 to 95–96 mm \cite{grasby_comparative_1996}. This inconsistent pattern indicates that the differences might have been done to local layout adjustments rather than any geometric plan to equalize apparent size for a viewer looking upward. Other proposed examples of unequal lettering such as the ``Priene Inscription of Alexander the Great” (334--323 BC, British Museum inv. 403 L.M.) as discussed by Pennethorne \cite{Pennethorne1878} and 
 the inscription at the top of the arch of Titus (81 CE, Rome) where one side has regular lettering and the other as unequal lettering, are equally underwhelming with no clear connection to optical corrections. Similarly, there are plenty of examples of non-uniform lettering at eye level such as the interior wall of the adyton at Didyma.

\begin{figure}[ht!]
 \centering
\includegraphics[width=0.9\linewidth]{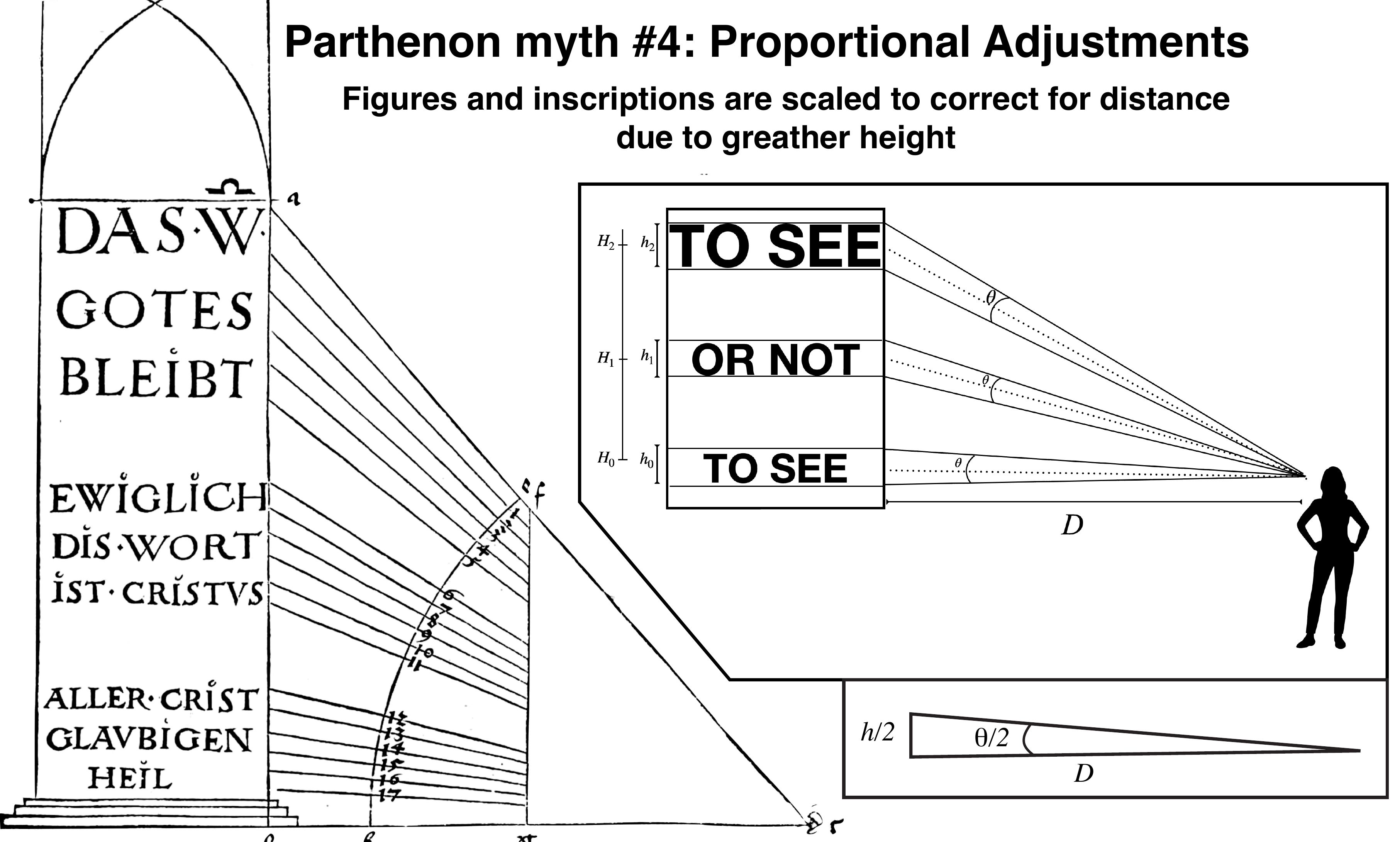}
 \caption{Correcting for height. If letters are scaled with height so that higher letters are larger than lower ones, they will have the same angular opening, hence project on the retina with the same size \cite{durer_underweysung_1538}.}
 \label{Fig-Durer}
\end{figure}

The scaling idea was revived and  widely discussed by many Renaissance authors concerned with optics and proportions. The initial example is found in  Albrecht Dürer's \textit{Underweysung der Messung} (1525) that provides a systematic geometric procedure for correction shown in  \cref{Fig-Durer}. His method, based on Euclidean optics, was subsequently illustrated and adapted by Joachim Fortius Ringelberg (1531), Sebastiano Serlio (1545), Walther Ryff (1547), Heinrich Lautensack (1564), Pietro Cataneo (1567), Daniele Barbaro (1569), Bartolomeo Romano (1595), and Gian Paolo Lomazzo (1584) \cite{Frangenberg1993}. All these authors have drawings similar to the one shown in \cref{Fig-Durer} and agreed that such corrections are necessary to preserve visual harmony and legibility in elevated architectural and sculptural settings.

\cor{In geometric terms, first place an object such as an inscribed character of height $h_0$ at eye level (see right-lower panel of \cref{Fig-Durer}), so that its midpoint is at the same height as the observing eye. Then the visual angle $\theta$  subtended by  this object placed at a distance $D$ is}
\begin{align}
\cor{    \theta = 2 \arctan \left( \frac{h_0}{2D} \right),}
\end{align}
\cor{Now, if we want to preserve the same angle $\theta$ (so that the retinal image has the same height as the initial object) when  an object is placed at a higher vertical position $H$, its height $h$ must increase accordingly.}
\cor{Explicitly, if the middle of a text is at height $H$ and the horizontal distance is $D$, then to keep $\theta$ constant, the height of this text must be
\begin{align}
 h=2 \cot(\theta)
   \left[\sqrt{\left(D^2+H^2\right)
   \tan ^2(\theta)
   +D^2}-D\right].
\end{align}}
\cor{In terms of $h_0$, the height is
\begin{align}
    h=\frac{{h_0}^2-4 D^2+\sqrt{16 {D}^4+8
   {D}^2 {h_0}^2+16 H^2
   {h_0}^2+{h_0}^4}}{2
   {h_0}}
\end{align}
We can now evaluate the correction needed by expanding this expression when $H/D$ is small (i.e. when the observer is sufficiently far compared to the height $H$):
\begin{align}
    h\approx {h_0}\left(1+\frac{4 H^2 }{4 {D}^2+{h_0}^2}+O\left((H/D)^3\right)\right).
\end{align}}
\cor{This last relation can be used to verify if such a correction is present from data by plotting $\log(h/h_0-1)$ as a function of $\log H$. On that graphs, point should be close to a line of slope $2$. The distance $D$ can then be determined from the intercept of this line with the vertical axis. No such scaling has ever been reported.}
 Importantly, since both $H$ and $h$ are fixed once engraved, there is one and  only one distance $D$ from which the angle $\theta$ is independent of $H$, hence only a single distance from which the effect of equal lettering is realized.

Because of this strong constraint about the vantage point,  these theoretical prescriptions promoted by so many authors were rarely, if ever, actually implemented in architectural, sculptural, or pictorial elements  as  Frangenberg argues \cite{Frangenberg1993}. The architect Scamozzi in \textit{L’Idea della Architettura Universale} (1615) explicitly dismissed Vitruvian-style adjustments as useless. Even in the Michelangelo’s painting, the \textit{Last Judgment},  found in the Palazzo Farnese where different rows of human figures have different sizes,  the evidence suggests intuitive modification rather than systematic application of a correction. Further, the higher row does appear to observers bigger from all viewpoints. The difficulty of finding any compelling implementation of this correction suggests that the notion of proportional adjustment  for height has remained largely a theoretical construct. This is mainly due to the fact that the human visual system naturally corrects for perspective effects. 

In psychology, \textit{size constancy} is the perceptual mechanism  by which the visual system maintains a stable sense of an object’s size despite changes in its retinal image caused by distance or perspective. As the distance $d$ to an object increases, the subtended angle $\theta$ decreases, meaning the retinal image becomes smaller. Yet in normal viewing we do not perceive distant objects as physically shrinking, because the brain integrates distance cues, such as convergence, parallax, texture gradients, shades, and linear perspective, to infer the object’s true scale. Therefore, when someone is on a balcony three floors up, they do not appear as giants when viewed from the ground despite the fact that their head is proportionally smaller on the retinal image. We naturally and subconsciously adjust our sense of proportions according to the angle of view and distances. Architects and artists who spend a long time looking at how objects and buildings are perceived in the real world are fully aware of perspective effects and see no reason to follow the advice of well-meaning theorists.

\section{Discussion}

\subsection{Fallacies}
All these myths rest on a number of fallacies that can be roughly classified as follows:\\

\noindent \textbf{The supreme knowledge fallacy.} The modern view of the Greek civilization as a pinnacle of knowledge and art took shape during the Enlightenment, when ancient Greece was celebrated as the birthplace of reason, democracy, and mathematics, and later intensified during the Romantic period, when the Greeks were idealized as embodying a lost unity of science, art, and philosophy \cite{Guthenke2008}.  To this day, intellectuals often project onto this ancient civilization an image of unrivaled and almost magical mastery across all domains of human knowledge, portraying them simultaneously as supreme mathematicians, engineers, psychologists, and artists whose collective genius culminated in monuments such as the Parthenon. Based on this romantic vision, the Parthenon is assumed to have been designed and executed with a perfect awareness of abstract geometry, structural mechanics, visual perception, and artistic theory, as though the builders possessed the combined knowledge of twentieth century modern sciences and more. Yet, no historical evidence supports these assumptions. As far as we know, the Greeks of the 5th century BCE  did not possess a theoretical science of engineering, or even a quantitative notion of force and buckling, nor did they have yet a mathematical formalization of conic sections, often  invoked in reconstructions of their design methods. More importantly, there is no evidence that they held any  knowledge of the particular visual illusions later alleged to justify the Parthenon’s refinements. 
Indeed, during the Parthenon’s building, the dominant theories of vision were those of natural philosophers such as Democritus, Plato, and Aristotle, who debated whether vision arose from emanations, internal fire, or the transparent properties of air and water \cite[p. 229]{VitruviusRowland1999}, hardly a basis to understand optical illusions. It is only a century later that Euclid’s Optics emerged with a more geometric view of vision. \\

\noindent \textbf{The intention fallacy.} The Parthenon myths rely on the implicit notion that we know the intent of its designers and what they tried to achieve.  Many authors have assumed that the intention behind certain curves is for the building to appear straight (\textit{linearity assumption}) or for certain irregularities in size to appear regular (\textit{regularity assumption}). Other authors have argued that the presence of curved elements in the building are meant to give an organic or corporeal feel to the building (\textit{the organic assumption}). Unfortunately, there is no contemporary record of any intentions behind the features of the Parthenon (or for that matter any classical Greek temples or constructions). We simply do not know the thought process of its creators or if they meant for the building to appear straight, curved,  perfectly regular,  slightly irregular, or evoke any particular feeling.\\

\noindent\textbf{The point(s)-of-view fallacy.} Many optical corrections, real or imaginary, require the observer to look at the building within a very specific region so that the projection of the temple on our retina has certain characteristics only visible from that limited vantage point \cite{MayerHillebrand1942,raymond1909proportion}. Away from this point, the explanation typically falls apart and becomes contradictory. For instance, the stylobate clearly looks curved from the corner. At what place should it look straight? Moreover, the location for a correction depends on the particular individual and different alleged corrections would require different points of view. For instance, the places poorly described by Pennethorne \cite{Pennethorne1878} where the Parthenon has perfect proportions is outside the purple areas found in \cref{corner-columns2}.  There is no evidence on the ground of the Acropolis that such special points ever existed (such as a marked point on the ground or a designated viewing platform, instructing observers to refrain from looking until reaching that exact position, lest they perceive the building imperfectly). It is very hard to imagine that any architect would build a temple with the express idea that its beauty can only be revealed at a few specific locations on the ground.\footnote{Note, however, that the architects of my own Institute, The Andrew Wiles Building in Oxford, played intentionally with columns of different thicknesses so that at a very particular vantage point, close to the Triton fountain,  they all align perfectly.} \\

\noindent \textbf{The illusion fallacy.} We know many different optical illusions. In particular some simple and very strong geometric illusions can be drawn on a piece of a paper.  These illusions have been mostly unveiled in the nineteenth century and have been since then studied extensively. The fallacy consists in believing that either some illusions exist while they have never been observed (such as a straight horizontal line appearing to sag under its own weight, or a straight column to appear hollow) or that these planar illusions are perceived the same way in the world around us and with the same magnitude. Indeed, it is well known that when we judge a three-dimensional scene, we integrate many different information to judge its physical and geometric attributes. We are not robots that extract outlines of an image and disregard all other relevant information about it. Moreover, our visual system is naturally tuned to three-dimensional scenes since it is the world in which we live. Often visual illusions on papers arise precisely because our cognitive system encodes assumptions about the outside world that cease to be valid in two dimensions (such as our best estimation of angles). It is therefore particularly ironic that multiple authors have used the argument that such illusions exist on paper to claim that they must exist in the real world, whereas they exist on paper precisely because they do not exist in the real world. \\

\noindent\textbf{The bootstrap fallacy.} An amusing argument for the existence of some optical corrections rests on the presumed existence of some other optical corrections, which allegdely proves that  the Greeks must have been aware of such an effect. Therefore, \textit{any} refinement can be explained by the same effect. If they corrected one illusion, surely they would have attempted to correct the other ones. The self-referential justification of the theory pulls itself by its own bootstraps. Once this argument is ignored and each optical correction is studied individually, the entire idea of optical corrections quickly falls apart. \\
  

\subsection{An alternative explanation} Demonstrating the absence of optical corrections does not explain the presence of architectural refinements. It is insufficient to dismiss them as arbitrary as our human thirst for meaning demands an explanation. A more pragmatic interpretation is that such refinements reflect empirical practices developed by builders and craftsmen. Throughout the history of architecture, empirical methods without explicit theoretical or scientific justification have been common across cultures. In China, the codified system of Feng Shui dictated the orientation and form of buildings according to geomantic ``energy flows'' that have no  physical basis \cite{lip1995,singh2019}. In Japan, hollow jars and charms were buried under temples and pagodas to ward off earthquakes. It worked with great success but only because the temple actual seismic resilience derived from joinery and the flexibility of timber frames \cite{coaldrake1996}. Medieval and Renaissance Europe preserved the Vitruvian idea of embedding resonant vessels in churches to improve acoustics, yet modern measurements demonstrate that their effect may be acoustically negligible \cite{godman2008enigma,valiere2013}. Similarly, Islamic architecture occasionally attributed acoustic powers to muqarnas vaulting \cite{necipoqlu1995}. Before Franklin’s experiments, European towers bore iron crosses or sacred talismans thought to repel lightning, though these provided no real protection \cite{lip1995}.  These examples illustrate that, across cultures, builders frequently relied on inherited empirical rules or symbolic reasoning that endured as part of architectural practice despite lacking scientific background. As the Greeks had not yet developed a theory of structures or mechanics, their architecture must have depended largely on such empirical conventions. It is therefore plausible that ancient Greek builders employed a repertoire of established empirical rules, attributed with specific purposes but devoid of theoretical foundation, the details of which are now sadly lost.

\subsection{Weighing the evidence}
We can now weigh the entirety of the evidence for the millennium-old myth of optical correction. First,  the observations and measurements clearly establish the existence of many refinements such as curved structures and irregular placements. Second, the entire historical evidence for optical corrections is based  on Vitruvius's writings, an author notoriously known to have misunderstood many aspects of Greek architecture \cite{Rankin1986}. It rests on a handful of opaque sentences written four centuries after the construction of the Parthenon, often interpreted and translated with a marked bias toward the theory of optical correction. Furthermore, modern perceptual science has shown that the alleged visual illusions either do not exist, only exist on paper, or have never been observed in situ.

Visual illusions are not a mysterious response that can only be appreciated through feelings, they are measurable physiological responses to optical stimuli, mostly due to signal processing by our visual system. In extracting information from pixels and orientation data that are triggered at the retina, the visual system makes some approximations that lead to under or over-evaluation of orientation and lengths. As such, if an effect exists, it is not a matter of opinion. Despite wide individual variations, illusions can be measured carefully in laboratories, or in situ, and can be estimated theoretically. Hence, the proponents of a theory of optical corrections should not turn to the prose of Vitruvius  but to experimental psychologists. If they believe an effect exists, then it needs to be properly documented and reproduced. The Parthenon copy at Nashville, Tennessee would be a good place to carry out these experiments. No anecdotal observations or contrived historical interpretations are sufficient to justify such an extraordinary claim.

\subsection{Conclusion}

The Parthenon is more than a building sitting on the Athenian Acropolis; it has become over the ages the very paradigm of antiquity, a vessel into which successive generations have poured their ideals of perfection. As such, it is the mother of all myths, each reflecting the intellectual preoccupations of the age. Vitruvius had promoted the idea of perfect proportional systems, an argument later dismantled by Claude Perrault \cite{Perrault1673}. In the nineteenth century, scholars and enthusiasts claimed to find the golden ratio within its proportions \cite{zeising1854lehre,ghyka1927esthetique}, a notion that has since been debunked \cite{markowsky1992misconceptions,livio2008golden} but will probably linger forever in the public mind. Similarly, the Parthenon has been cast as an ultimate trompe-l'oeil, supposedly designed to exploit visual illusions in its curves and refinements. Yet no historical or scientific evidence support such claims. The main lesson from the Parthenon is that we are prone to delude ourselves in the pursuit of ancient mysteries and ideal fantasies of order and harmony.\\

There is, nevertheless, much to be learned from the wisdom of ancient Greece. In the \textit{Republic} (Book VII, 523a–527c), Plato insists that geometry trains the soul to see beyond sensory deception. It is time to transcend our own tendency to see illusions where none exist and apply mathematics and science to see beyond our own delusions.\\




\noindent\textbf{Acknowledgments:} It is my great pleasure to acknowledge the extraordinary help of Dr Nicola Kirkham in gathering key ancient texts and extracting data from them.

 \bibliographystyle{RS} 

\end{document}